\documentclass[11pt]{amsart}
\usepackage{amssymb,mathrsfs}
\title[\tiny {The Schr{\"o}dinger-Weil Representation and Jacobi Forms of Half Integral Weight }]
{The Schr{\"o}dinger-Weil Representation and Jacobi Forms of Half Integral Weight }

\author{Jae-Hyun Yang}
\address{Department of Mathematics, Inha University, Incheon
402-751, Korea}
\email{jhyang@inha.ac.kr }

\begin{document}

\newtheorem{theorem}{Theorem}[section]
\newtheorem{corollary}{Corollary}[section]
\newtheorem{lemma}{Lemma}[section]
\newtheorem{proposition}{Proposition}[section]
\newtheorem{remark}{Remark}[section]
\newtheorem{definition}{Definition}[section]

\renewcommand{\theequation}{\thesection.\arabic{equation}}
\renewcommand{\thetheorem}{\thesection.\arabic{theorem}}
\renewcommand{\thelemma}{\thesection.\arabic{lemma}}
\newcommand{\BR}{\mathbb R}
\newcommand{\BQ}{\mathbb Q}
\newcommand{\bn}{\bf n}
\def\charf {\mbox{{\text 1}\kern-.24em {\text l}}}
\newcommand{\BC}{\mathbb C}
\newcommand{\BZ}{\mathbb Z}

\thanks{\noindent{2010 Mathematics Subject Classification:} Primary 11F27, 11F50\\
\indent  Keywords and phrases: the Schr{\"o}dinger-Weil
representation, covariant maps, Jacobi forms}

%\maketitle

\begin{abstract}
{In this paper, we define the concept of Jacobi forms of half integral weight using Takase's automorphic factor
of weight $1/2$ for a two-fold covering group of the symplectic group $Sp(n,{\mathbb R})$ on the Siegel
upper half plane and find covariant maps for the Schr{\"o}dinger-Weil representation. Using these covariant maps, we construct
Jacobi forms of half integral weight with respect to an arithmetic subgroup of the Jacobi group.}
\end{abstract}
\maketitle

%\tableofcontents
\newcommand\tr{\triangleright}
\newcommand\al{\alpha}
\newcommand\be{\beta}
\newcommand\g{\gamma}
\newcommand\gh{\Cal G^J}
\newcommand\G{\Gamma}
\newcommand\de{\delta}
\newcommand\e{\epsilon}
\newcommand\z{\zeta}
\newcommand\vth{\vartheta}
\newcommand\vp{\varphi}
\newcommand\om{\omega}
\newcommand\p{\pi}
\newcommand\la{\lambda}
\newcommand\lb{\lbrace}
\newcommand\lk{\lbrack}
\newcommand\rb{\rbrace}
\newcommand\rk{\rbrack}
\newcommand\s{\sigma}
\newcommand\w{\wedge}
\newcommand\fgj{{\frak g}^J}
\newcommand\lrt{\longrightarrow}
\newcommand\lmt{\longmapsto}
\newcommand\lmk{(\lambda,\mu,\kappa)}
\newcommand\Om{\Omega}
\newcommand\ka{\kappa}
\newcommand\ba{\backslash}
\newcommand\ph{\phi}
\newcommand\M{{\Cal M}}
\newcommand\bA{\bold A}
\newcommand\bH{\bold H}

\newcommand\Hom{\text{Hom}}
\newcommand\cP{\Cal P}
\newcommand\cH{\Cal H}

\newcommand\pa{\partial}

\newcommand\pis{\pi i \sigma}
\newcommand\sd{\,\,{\vartriangleright}\kern -1.0ex{<}\,}
\newcommand\wt{\widetilde}
\newcommand\fg{\frak g}
\newcommand\fk{\frak k}
\newcommand\fp{\frak p}
\newcommand\fs{\frak s}
\newcommand\fh{\frak h}
\newcommand\Cal{\mathcal}

\newcommand\fn{{\frak n}}
\newcommand\fa{{\frak a}}
\newcommand\fm{{\frak m}}
\newcommand\fq{{\frak q}}
\newcommand\CP{{\mathcal P}_g}
\newcommand\Hgh{{\mathbb H}_g \times {\mathbb C}^{(h,g)}}
\newcommand\BD{\mathbb D}
\newcommand\BH{\mathbb H}
\newcommand\CCF{{\mathcal F}_g}
\newcommand\CM{{\mathcal M}}
\newcommand\Ggh{\Gamma_{g,h}}
\newcommand\Chg{{\mathbb C}^{(h,g)}}
\newcommand\Yd{{{\partial}\over {\partial Y}}}
\newcommand\Vd{{{\partial}\over {\partial V}}}

\newcommand\Ys{Y^{\ast}}
\newcommand\Vs{V^{\ast}}
\newcommand\LO{L_{\Omega}}
\newcommand\fac{{\frak a}_{\mathbb C}^{\ast}}

\renewcommand\th{\theta}
\renewcommand\l{\lambda}
\renewcommand\k{\kappa}
\newcommand\tg{\widetilde\gamma}
\newcommand\wmo{{\mathscr W}_{\mathcal M,\Omega}}
\newcommand\hrnm{H_\BR^{(n,m)}}
\newcommand\rmn{\BR^{(m,n)}}

\begin{section}{{\bf Introduction}}
\setcounter{equation}{0}

For a given fixed positive integer $n$, we let
$${\mathbb H}_n=\,\big\{\,\Om\in \BC^{(n,n)}\,\big|\ \Om=\,^t\Om,\ \ \ \text{Im}\,\Om>0\,\big\}$$
be the Siegel upper half plane of degree $n$ and let
$$Sp(n,\BR)=\big\{ g\in \BR^{(2n,2n)}\ \big| \ ^t\!gJ_ng= J_n\ \big\}$$
be the symplectic group of degree $n$, where $F^{(k,l)}$ denotes
the set of all $k\times l$ matrices with entries in a commutative
ring $F$ for two positive integers $k$ and $l$, $^t\!M$ denotes
the transpose of a matrix $M,\ \text{Im}\,\Om$ denotes the
imaginary part of $\Om$ and
$$J_n=\begin{pmatrix} 0&I_n \\
                   -I_n&0 \\ \end{pmatrix}.$$
Here $I_n$ denotes the identity matrix of degree $n$.
We see that $Sp(n,\BR)$ acts on $\BH_n$ transitively by
\begin{equation*}g\cdot \Om=(A\Om+B)(C\Om+D)^{-1}, \end{equation*}
where $g=\begin{pmatrix} A&B\\ C&D\end{pmatrix}\in Sp(n,\BR)$ and
$\Om\in \BH_n.$

For two positive integers $n$ and $m$, we consider the Heisenberg
group
$$H_{\BR}^{(n,m)}=\{\,(\l,\mu;\k)\,|\ \l,\mu\in \BR^{(m,n)},\ \k\in \BR^{(m,m)},\ \
\k+\mu\,^t\l\ \text{symmetric}\ \}$$ endowed with the following
multiplication law
$$(\l,\mu;\k)\circ (\l',\mu';\k')=(\l+\l',\mu+\mu';\k+\k'+\l\,^t\mu'-
\mu\,^t\l').$$ We let
$$G^J=Sp(n,\BR)\ltimes H_{\BR}^{(n,m)}\quad \ ( \textrm{semi-direct product})$$
be the Jacobi group endowed with the following multiplication law
$$\Big(g,(\lambda,\mu;\kappa)\Big)\cdot\Big(g',(\lambda',\mu';\kappa')\Big) =\,
\Big(gg',(\widetilde{\lambda}+\lambda',\widetilde{\mu}+ \mu';
\kappa+\kappa'+\widetilde{\lambda}\,^t\!\mu'
-\widetilde{\mu}\,^t\!\lambda')\Big)$$ with $g,g'\in Sp(n,\BR),
(\lambda,\mu;\kappa),\,(\lambda',\mu';\kappa') \in
H_{\BR}^{(n,m)}$ and
$(\widetilde{\lambda},\widetilde{\mu})=(\lambda,\mu)g'$. We let
$\G_n=Sp(n,\BZ)$ be the Siegel modular group of degree $n$. We let
$$\G^J_n=\G_n\ltimes H_{\BZ}^{(n,m)}$$
be the Jacobi modular group. Then we have the {\it natural action}
of $G^J$ on the Siegel-Jacobi space $\BH_{n,m}:=\BH_n\times
\BC^{(m,n)}$ defined by
\begin{equation}\Big(g,(\lambda,\mu;\kappa)\Big)\cdot (\Om,Z)=\Big(g\!\cdot\! \Om,(Z+\lambda \,\Om+\mu)
(C\,\Om+D)^{-1}\Big),
\end{equation}
where $g=\begin{pmatrix} A&B\\
C&D\end{pmatrix} \in Sp(n,\BR),\ (\lambda,\mu; \kappa)\in
H_{\BR}^{(n,m)}$ and $(\Om,Z)\in \BH_{n,m}.$ We refer to \cite{BR, CR},\,\cite{EZ},\,
\cite{YJ6}-\cite{YJ15} for more details on materials, e.g., Jacobi forms, invariant metrics,
invariant differential operators, Maass-Jacobi forms etc related to
the Siegel-Jacobi space.

\vskip 0.2cm The Weil representation for a symplectic group was
first introduced by A. Weil in \cite{W} to reformulate Siegel's
analytic theory of quadratic forms (cf.\,\cite{Si}) in terms of
the group theoretical theory. It is well known that the Weil
representation plays a central role in the study of the
transformation behaviors of theta series.
Whenever we study the transformation formulas of theta series or Siegel modular forms of
half integral weights, we are troubled by the ambiguity of the factor $\det (C\Om+D)^{1/2}$ in its signature.
This means that we should consider the transformation formula on a non-trivial two-fold covering group of
a symplectic group. In his paper \cite{Ta3}, Takase removed the ambiguity of the factor $\det (C\Om+D)^{1/2}$
by constructing the {\it right\ explicit} automorphic factor $J_{1/2}$ of weight $1/2$ for $Sp(n,\BR)_*$ on
$\BH_n$\,:
\begin{equation}
J_{1/2}: Sp(n,\BR)_*\times \BH_n\lrt \BC^*.
\end{equation}
Here $Sp(n,\BR)_*$ is the two-fold covering group of $Sp(n,\BR)$ in the sense of a real Lie group.
See (4.15) for the precise definition. $J_{1/2}$ is real analytic on $Sp(n,\BR)_*$, holomorphic on
$\BH_n$ and satisfies the relation
\begin{equation}
J_{1/2}(g_*,\Om)^2=\det (C\Om+D),
\end{equation}
where $g_*=(g,\epsilon)\in Sp(n,\BR)_*$ with $g=\begin{pmatrix} A & B \\ C & D \end{pmatrix}\in Sp(n,\BR).$
Using the automorphic factor $J_{1/2}$ of weight $1/2$, Takase expressed the transformation formula of theta series without
ambiguity of $\det (C\Om+D)^{1/2}$. Moreover he decomposed the automorphic factor
$$j(\gamma,\Om)= {{\vartheta(\g\cdot\Om)}\over {\vartheta(\Om)}}$$
with a standard theta series $\vartheta (\Om)$ into a product of a character and the automorphic factor
$J_{1/2}(\g_*,\Om)$. The automorphic factor $J_{1/2}$ of Takase will play an important role in the further study of
half integral weight Siegel modular forms and half integral weight Jacobi forms.

\vskip 0.21cm
This paper is organized as follows.
In Section 2, we discuss the Schr{\"o}dinger representation of the Heisenberg group
$H_\BR^{(n,m)}$ associated with a nonzero symmetric real matrix of
degree $m$ which is formulated in \cite{YJ1, YJ2}.
In Section 3,
we define the Schr{\"o}dinger-Weil
representation $\omega_\CM$ of the Jacobi group $G^J$ associated
with a symmetric positive definite matrix $\CM$ and provide some
of the actions of $\omega_\CM$ on the representation space
$L^2\big(\BR^{(m,n)}\big)$ explicitly.
In Section 4, we review Jacobi forms of integral weight, Siegel modular forms of half integral weight,
and define Jacobi forms of half integral weight using the automorphic factor $J_{1/2}$ of weight $1.2$ for the
metaplectic group $Sp(n,\BR)_*$ on the Siegel upper half plane.
In Section 5, we find covariant maps for the Schr{\"o}dinger-Weil
representation $\omega_\CM$.
In the final section we construct Jacobi forms of half integral weight with
respect to an arithmetic subgroup of $\Gamma^J$ using
covariant maps obtained in Section 5.

\vskip 0.2cm \noindent {\bf Notations\,:} \ \ We denote by $\BZ,\,\,\BR$
and $\BC$ the ring of integers, the field of real numbers and the field of complex numbers
respectively. $\BC^{\times}$ denotes the multiplicative group of
nonzero complex numbers and $\BZ^{\times}$ denotes the set of all nonzero integers..
$T$ denotes the multiplicative group of
complex numbers of modulus one. The symbol ``:='' means that the
expression on the right is the definition of that on the left. For
two positive integers $k$ and $l$, $F^{(k,l)}$ denotes the set of
all $k\times l$ matrices with entries in a commutative ring $F$.
For a square matrix $A\in F^{(k,k)}$ of degree $k$, $\sigma(A)$
denotes the trace of $A$. For any $M\in F^{(k,l)},\ {}^t\!M$ denotes
the transpose of a matrix $M$. $I_n$ denotes the identity matrix of
degree $n$. We put $i=\sqrt{-1}.$ For $z\in\BC,$ we define
$z^{1/2}=\sqrt{z}$ so that $-\pi / 2 < \ \arg (z^{1/2})\leqq
\pi/2.$ Further we put $z^{\kappa/2}=\big(z^{1/2}\big)^\kappa$ for
every $\kappa\in\BZ.$
For a positive integer $m$ we denote by $S(m,F)$ the additive group consisting of all $m\times m$
symmetric matrices with coefficients in a commutative ring $F$.

\end{section}

\vskip 1cm

\begin{section}{{\bf The Schr{\"o}dinger Representation}}
% of $H_\BR^{(n,m)}$}}
\setcounter{equation}{0}

\vskip 0.2cm
%\vspace{0.1in}\\
%\indent
First of all, we observe that $H_{\mathbb{R}}^{(n,m)}$ is a 2-step
nilpotent Lie group. The inverse of an element
$(\lambda,\mu;\kappa)\in H_{\mathbb {R}}^{(n,m)}$ is given by
$$(\lambda,\mu;\kappa)^{-1}=(-\lambda,-\mu;-\kappa+\lambda\,^t\!\mu-\mu\,^t\!\lambda).$$
Now we set
\begin{equation*}
[\lambda,\mu;\kappa]=(0,\mu;\kappa)\circ
(\lambda,0;0)=(\lambda,\mu;\kappa-\mu\,^t\! \lambda).
\end{equation*}
\noindent Then $H_{\mathbb {R}}^{(n,m)}$ may be regarded as a
group equipped with the following multiplication
\begin{equation*}
[\lambda,\mu;\kappa]\diamond
[\lambda_0,\mu_0;\kappa_0]=[\lambda+\lambda_0,\mu+\mu_0;
\kappa+\kappa_0+\lambda\,^t\!\mu_0+\mu_0\,^t\!\lambda].
\end{equation*}
\noindent The inverse of $[\lambda,\mu;\kappa]\in H_{\mathbb
{R}}^{(n,m)}$ is given by
$$[\lambda,\mu;\kappa]^{-1}=[-\lambda,-\mu;-\kappa+\lambda\,^t\!\mu+\mu\,^t\!\lambda].$$
We set
\begin{equation*}
 L=\left\{\,[0,\mu;\kappa]\in H_{\mathbb
{R}}^{(n,m)}\,\Big| \, \mu\in \mathbb {R}^{(m,n)},\
\kappa=\,^t\!\kappa\in \mathbb {R}^{(m,m)}\ \right\}.
 \end{equation*}
\noindent Then $L$ is a commutative normal subgroup of $H_{\mathbb
{R}}^{(n,m)}$. Let ${\widehat {L}}$ be the Pontrajagin dual of
$L$, i.e., the commutative group consisting of all unitary
characters of $L$. Then ${\widehat {L}}$ is isomorphic to the
additive group $\mathbb {R}^{(m,n)}\times S(m,\mathbb
{R})$ via
\begin{equation*}
\langle a,{\hat a}\rangle =e^{2 \pi i\sigma({\hat
{\mu}}\,^t\!\mu+{\hat {\kappa}}\kappa)}, \ \ \ a=[0,\mu;\kappa]\in
L,\ {\hat {a}}=({\hat {\mu}},{\hat {\kappa}})\in \widehat{L},
\end{equation*}

\noindent where $S(m,\mathbb {R})$ denotes the space of
all symmetric $m\times m$ real matrices.

\vskip 0.21cm
We put
\begin{equation*}
S=\left\{\,[\lambda,0;0]\in H_{\mathbb {R}}^{(n,m)}\,\Big|\
\lambda\in \mathbb {R}^{(m,n)}\, \right\}\cong \mathbb
{R}^{(m,n)}.
\end{equation*}
\noindent Then $S$ acts on $L$ as follows:
\begin{equation*}
\alpha_{\lambda}([0,\mu;\kappa])=[0,\mu;\kappa+\lambda\,^t\!\mu+\mu\,^t\!\lambda],
\ \ \ [\lambda,0,0]\in S.
\end{equation*}
\noindent We see that the Heisenberg group $\left( H_{\mathbb
{R}}^{(n,m)}, \diamond\right)$ is isomorphic to the semi-direct
product $S\ltimes L$ of $S$ and $L$ whose multiplication is given
by
$$(\lambda,a)\cdot
(\lambda_0,a_0)=\big(\lambda+\lambda_0,a+\alpha_{\lambda}(a_0)\big),\
\ \lambda,\lambda_0\in S,\ a,a_0\in L.$$
On the other hand, $S$
acts on ${\widehat {L}}$ by
\begin{equation*}
\alpha_{\lambda}^{*}({\widehat {a}})=({\hat {\mu}}+2{\hat
{\kappa}}\lambda, {\hat {\kappa}}),\ \ [\lambda,0;0]\in S,\ \
{\widehat a}=({\hat {\mu}},{\hat {\kappa}})\in {\widehat {L}}.
\end{equation*}
\noindent Then, we have the relation $\langle
\alpha_{\lambda}(a),{\widehat {a}} \rangle=\langle
a,\alpha_{\lambda}^{*} ({\widehat {a}}) \rangle$ for all $a\in L$ and
${\widehat {a}}\in {\widehat {L}}.$

\vskip 0.35cm
 \indent We have three types of $S$-orbits in ${\widehat
{L}}.$

\smallskip

\noindent {\scshape Type I.} Let ${\hat{\kappa}} \in
S(m,\mathbb {R})$ be nondegenerate. The $S$-orbit of
${\hat {a}}({\hat {\kappa}})=(0,{\hat {\kappa}}) \in {\widehat
{L}}$ is given by

\begin{equation*}
\widehat{\mathcal{O}}_{\hat{\kappa}}= \left\{(2{\hat
{\kappa}}\lambda,{\hat {\kappa}}) \in {\widehat {L}}\ \Big|\
\lambda \in \mathbb{R}^{(m,n)} \right\} \cong \mathbb {R}^{(m,n)}.
\end{equation*}

\noindent {\scshape Type II.}\ \ Let
$({\hat{\mu}},{\hat{\kappa}})\in
\mathbb{R}^{(m,n)}\times  S(m,\mathbb {R})$ with ${\hat \mu}\in \BR^{(m,n)},\
{\hat\kappa}\in S(m,\BR)$ and
degenerate  ${\hat{\kappa}}\neq 0.$ Then
\begin{equation*}
\widehat{\mathcal{O}}_{(\hat{\mu}, \hat{\kappa})} = \left\{ (\hat
\mu + 2\hat{\kappa}\lambda, \hat{\kappa}) \Big|\ \lambda \in
\mathbb {R}^{(m,n)} \right\} \subsetneqq \mathbb {R}^{(m,n)}\times
\{ \hat{\kappa} \} .
\end{equation*}

\noindent {\scshape Type III.} Let $\hat{y} \in
\mathbb{R}^{(m,n)}$. The $S$-orbit ${\widehat {\mathcal {O}
}}_{\hat {y}}$ of $\hat{a}(\hat{y}) = (\hat{y} ,0)$ is given by
\begin{equation*}
{\widehat {\mathcal {O} }}_{\hat {y}}=\left\{\,({\hat
{y}},0)\,\right\}={\hat {a}} ({\hat {y}}).
\end{equation*}
\noindent We have
$$
{\widehat{L}}= \left( \bigcup_{\begin{subarray}{c} \hat{\kappa}
\in S(m,\mathbb{R}) \\ {\hat{\kappa}} \,\
\text{nondegenerate} \end{subarray}}
 \widehat {\mathcal{O}}_ {\hat\kappa }
 \right)
  \bigcup
\left( \bigcup_{{\hat {y}}\in \mathbb{R}^{(m,n)}}{\widehat
{\mathcal {O} }}_{\hat {y}}\right)  \bigcup  \left(
\bigcup_{\begin{subarray}{c}({\hat{\mu}},{\hat {\kappa}}) \in
\mathbb{R}^{(m,n)} \times S(m,\mathbb{R}) \\
\hat{\kappa} \neq 0 \,\ \text{degenerate} \end{subarray}}
 {\widehat {\mathcal {O}}}_{({\hat \mu},{\hat \kappa })} \right)
$$
\noindent as a set. The stabilizer $S_{\hat {\kappa}}$ of $S$ at
${\hat {a}}({\hat {\kappa}})=(0,{\hat {\kappa}})$ is given by
\begin{equation*}
S_{\hat {\kappa}}=\{0\}.
 \end{equation*}
\noindent And the stabilizer $S_{\hat {y}}$ of $S$ at ${\hat
{a}}({\hat {y}})= ({\hat {y}},0)$ is given by
\begin{equation*}
S_{\hat {y}}=\left\{\,[\lambda,0;0]\,\Big|\ \lambda\in \mathbb
{R}^{(m,n)}\,\right\}=S \,\cong\,\mathbb {R}^{(m,n)}.
\end{equation*}
\indent In this section, for the present being we set
$H=H_{\mathbb {R}}^{(n,m)}$ for brevity. We see that $L$ is a
closed, commutative normal subgroup of $H$. Since
$(\lambda,\mu;\kappa)=(0,\mu; \kappa+\mu\,^t\!\lambda)\circ
(\lambda,0;0)$ for $(\lambda,\mu;\kappa)\in H,$ the homogeneous
space $X=L\backslash H$ can be identified with
$\mathbb{R}^{(m,n)}$ via
$$
Lh=L\circ (\lambda,0;0)\longmapsto \lambda,\ \ \
h=(\lambda,\mu;\kappa)\in H.
$$
We observe that $H$ acts on $X$ by
\begin{equation*}
(Lh)\cdot h_0=L\,(\lambda+\lambda_0,0;0)=\lambda+\lambda_0,
\end{equation*}
\noindent where $h=(\lambda,\mu;\kappa)\in H$ and
$h_0=(\lambda_0,\mu_0;\kappa_0)\in H.$

\medskip

\indent If $h=(\lambda,\mu;\kappa)\in H$, we have
\begin{equation*}
l_h=(0, \mu; \kappa+\mu\,^t\!\lambda),\ \ \ s_h=(\lambda,0;0)
\end{equation*}
\noindent in the Mackey decomposition of $h=l_h \circ s_h$
(cf.\,\cite{M}). Thus if $h_0=(\lambda_0,\mu_0;\kappa_0)\in H,$
then we have
\begin{equation*}
s_h\circ h_0=(\lambda,0;0)\circ
(\lambda_0,\mu_0;\kappa_0)=(\lambda+\lambda_0,\mu_0;
\kappa_0+\lambda\,^t\!\mu_0)
\end{equation*}
\noindent and so
\begin{equation}
l_{s_h\circ
h_0}=\big(0,\mu_0;\kappa_0+\mu_0\,^t\!\lambda_0+\lambda\,^t\!\mu_0+
\mu_0\,^t\!\lambda\big).
\end{equation}
 \indent For a real symmetric matrix
$c=\,^tc\in S(m,\BR) $ with $c\neq 0$, we consider the
unitary character $\chi_c$ of $L$ defined by
\begin{equation}
\chi_c\left((0,\mu;\kappa)\right)=e^{\pi i \sigma(c\kappa)}\,I,\ \
\ (0,\mu;\kappa)\in L,
\end{equation}
\noindent where $I$ denotes the identity mapping. Then the
representation ${\mathscr W}_c=\text{Ind}_L^H\,\chi_c$ of $H$
induced from $\chi_c$ is realized on the Hilbert space
$H(\chi_c)=L^2\big(X,d{\dot {h}},\mathbb {C}\big) \cong
L^2\left(\mathbb{R}^{(m,n)}, d\xi\right)$ as follows. If
$h_0=(\lambda_0,\mu_0;\kappa_0)\in H$ and $x=Lh\in X$ with
$h=(\lambda,\mu;\kappa)\in H,$ we have
\begin{equation}
\left({\mathscr W}_c (h_0)f\right)(x)=\chi_c (l_{s_h\circ
h_0}) f(xh_0),\ \ f\in H(\chi_c).
\end{equation}
\noindent It follows from (2.1) that
\begin{equation}
\left( {\mathscr W}_c (h_0)f\right)(\lambda)=e^{\pi
i\sigma\{c(\kappa_0+\mu_0\,^t\!\lambda_0+
2\lambda\,^t\!\mu_0)\}}\,f(\lambda+\lambda_0),
\end{equation}

\noindent where $h_0=(\lambda_0,\mu_0;\kappa_0)\in H$ and
$\lambda\in\BR^{(m,n)}.$ Here we identified
$x=Lh$\,(resp.\,$xh_0=Lhh_0$) with $\lambda$\,(resp.\,
$\lambda+\lambda_0$). The induced representation ${\mathscr W}_c$
is called the $\textsf{Schr{\" {o}}dinger\ representation}$ of $H$
associated with $\chi_c.$ Thus ${\mathscr W}_c $ is a monomial
representation.

\medskip

\begin{theorem}
 Let $c$ be a positive definite symmetric real matrix of
degree $m$. Then the Schr{\" {o}}dinger representation ${\mathscr
W}_c $ of $H$ is irreducible. \end{theorem} \noindent{\it Proof.}\
 The proof can be found in \cite{YJ1},\ Theorem 3. \hfill$\square$

\vskip 0.2cm\noindent {\bf Remark.} We refer to
\cite{YJ1}-\cite{YJ5} for more representations of the Heisenberg
group $H_\BR^{(n,m)}$ and their related topics.

\end{section}

\vskip 1cm

\begin{section}{{\bf The Schr{\"o}dinger-Weil Representation}}
\setcounter{equation}{0}

\vskip 0.2cm Throughout this section we assume that $\CM$ is a
symmetric real positive definite $m\times m$ matrix. We
consider the Schr{\"o}dinger representation ${\mathscr W}_\CM$ of
the Heisenberg group $\hrnm$ with the central character ${\mathscr
W}_\CM((0,0;\kappa))=\chi_\CM((0,0;\kappa))=e^{\pi
i\,\s(\CM\kappa)},\ \kappa\in S(m,\BR)$\,(cf.\,(2.2)). We
note that the symplectic group $Sp(n,\BR)$ acts on $\hrnm$ by
conjugation inside $G^J$. For a fixed element $g\in Sp(n,\BR)$,
the irreducible unitary representation ${\mathscr W}_\CM^g$ of
$\hrnm$ defined by
\begin{equation}
{\mathscr W}_\CM^g(h)={\mathscr W}_\CM(ghg^{-1}),\quad h\in\hrnm
\end{equation}
has the property that

\begin{equation*}
{\mathscr W}_\CM^g((0,0;\k))={\mathscr W}_\CM((0,0;\k))=e^{\pi
i\,\s(\CM \k)}\,\textrm{Id}_{H(\chi_\CM)},\quad \k\in
S(m,\BR).
\end{equation*}
Here $\textrm{Id}_{H(\chi_\CM)}$ denotes the identity operator on
the Hilbert space $H(\chi_\CM).$ According to Stone-von Neumann
theorem, there exists a unitary operator $R_\CM(g)$ on
$H(\chi_\CM)$ such that
\begin{equation}
R_\CM(g){\mathscr W}_\CM(h)={\mathscr
W}_\CM^g(h) R_\CM(g)\qquad {\rm for\ all}\ h\in\hrnm.
\end{equation}
We observe that
$R_\CM(g)$ is determined uniquely up to a scalar of modulus one.

\vskip 0.35cm
From now on, for brevity, we put $G=Sp(n,\BR).$ According to
Schur's lemma, we have a map $c_\CM:G\times G\lrt T$ satisfying
the relation
\begin{equation}
R_\CM(g_1g_2)=c_\CM(g_1,g_2)R_\CM(g_1)R_\CM(g_2)\quad \textrm{for
all }\ g_1,g_2\in G.
\end{equation}
Therefore $R_\CM$ is a projective representation of $G$ on
$H(\chi_\CM)$ and $c_\CM$ defines the cocycle class in $H^2(G,T).$
The cocycle $c_\CM$ yields the central extension $G_\CM$ of $G$ by
$T$. The group $G_\CM$ is a set $G\times T$ equipped with the
following multiplication

\begin{equation}
(g_1,t_1)\cdot (g_2,t_2)=\big(g_1g_2,t_1t_2\,
c_\CM(g_1,g_2)^{-1/m}\,\big),\quad g_1,g_2\in G,\ t_1,t_2\in T.
\end{equation}
We see immediately that the map ${\widetilde R}_\CM:G_\CM\lrt
GL(H(\chi_\CM))$ defined by

\begin{equation}
{\widetilde R}_\CM(g,t)=t^m\,R_\CM(g) \quad \textrm{for all}\
(g,t)\in G_\CM
\end{equation}
is a true representation of $G_\CM.$ As in Section 1.7 in
\cite{LV}, we can define the map $s_\CM:G\lrt T$ satisfying the
relation

\begin{equation*}
c_\CM(g_1,g_2)^2=s_\CM(g_1)^{-1}s_\CM(g_2)^{-1}s_\CM(g_1g_2)\quad
\textrm{for all}\ g_1,g_2\in G.
\end{equation*}
Thus we see that

\begin{equation}
G_{2,\CM}=\left\{\, (g,t)\in G_\CM\,|\ t^2=s_\CM(g)^{-1/m}\,\right\}
\end{equation}

\noindent is the metaplectic group associated with $\CM$ that is a
two-fold covering group of $G$. The restriction $R_{2,\CM}$ of
${\widetilde R}_\CM$ to $G_{2,\CM}$ is the $\textsf{Weil representation}$ of
$G$ associated with $\CM$. Now we define the projective
representation $\pi_\CM$ of the Jacobi group $G^J$ by

\begin{equation}
\pi_\CM(hg)={\mathscr W}_\CM(h)\,R_\CM(g),\quad h\in\hrnm,\ g\in
G.
\end{equation}
Here we identified $h=(\lambda,\mu;\kappa)\in \hrnm$ (resp. $g\in Sp(n,\BR)$) with
$(I_{2n},(\lambda,\mu;\kappa))\in G^J$ (resp. $(g,(0,0;0))\in G^J).$
The projective representation $\pi_\CM$ of $G^J$ is naturally
extended to the true representation $\omega_\CM$ of the group
$G_{2,\CM}^J\!=G_{2,\CM}\ltimes \hrnm.$ The representation
$\om_\CM$ is called the $\textsf{Schr{\"o}dinger-Weil}$
$\textsf{representation}$ of $G^J.$ Indeed we have
\begin{equation}
\omega_\CM(h\!\cdot\!(g,t))=t^m\,\pi_\CM(hg)=t^m\, {\mathscr
W}_\CM(h)\,R_\CM(g),\quad h\in\hrnm,\ (g,t)\in G_{2,\CM}.
\end{equation}
Here we identified $h=(\lambda,\mu;\kappa)\in \hrnm$ (resp. $(g,t)\in G_{2,\CM})$ with
$((I_{2n},1),(\lambda,\mu;\kappa))\in G^J_{2,\CM}$ (resp. $((g,t),(0,0;0))\in G^J_{2,\CM}).$
\vskip0.2cm We recall that the following matrices
\begin{eqnarray*}
t(b)&=&\begin{pmatrix} I_n& b\\
                   0& I_n\end{pmatrix}\ \textrm{with any}\
                   b=\,{}^tb\in \BR^{(n,n)},\\
g(\alpha)&=&\begin{pmatrix} {}^t\alpha & 0\\
                   0& \alpha^{-1}  \end{pmatrix}\ \textrm{with
                   any}\ \alpha\in GL(n,\BR),\\
\s_n&=&\begin{pmatrix} 0& -I_n\\
                   I_n&\ 0\end{pmatrix}
\end{eqnarray*}
\noindent generate the symplectic group $G=Sp(n,\BR)$
(cf.\,\cite[p.\,326]{F},\,\cite[p.\,210]{Mum}). Therefore the
following elements $h_t(\lambda,\mu;\kappa),\
t_\CM(b;t),\,g_\CM(\alpha;t)$ and $\s_{n,\CM;t}$ of $G_\CM\ltimes \hrnm$
defined by
\begin{eqnarray*}
&& h_t(\la,\mu;\kappa)=\big( (I_{2n},t),(\la,\mu;\kappa)\big)\
\textrm{with}\ t\in T,\ \la,\mu\in
\BR^{(m,n)}\ \textrm{and}\ \kappa\in\BR^{(m,m)} ,\\
&&t_\CM(b;t)=\big( (t(b),t),(0,0;0) \big)\ \textrm{with any}\
                   b=\,{}^tb\in \BR^{(n,n)},\ t\in T,\\
&& g_\CM(\alpha;t)=\left(
\big(g(\alpha),t),(0,0;0)\right)\
\textrm{with any}\ \alpha\in GL(n,\BR)\ \textrm{and}\ t\in T,\\
 &&\s_{n,\CM;t}=\left( (\s_n,t),(0,0;0)\right)\
 \textrm{with}\ t\in T
\end{eqnarray*}
generate the group $G_\CM\ltimes\hrnm.$ We can show that the
representation ${\widetilde R}_\CM$ is realized on the
representation $H(\chi_\CM)=L^2\big(\rmn\big)$ as follows: for
each $f\in L^2\big(\rmn\big)$ and $x\in \rmn,$ the actions of
${\widetilde R}_\CM$ on the generators are given by

\begin{eqnarray}
%\begin{equation}
\left( {\widetilde R}_\CM
\big(h_t(\lambda,\mu;\kappa)\big)f\right)(x)&=&\,t^m\,e^{\pi
i\,\s\{\CM(\kappa+\mu\,{}^t\!\lambda+2\,x\,{}^t\mu)\}}\,f(x+\lambda),\\
%\end{equation}
%\begin{eqnarray}
\left( {\widetilde R}_\CM\big(t_\CM(b;t)\big)f\right)(x)&=& t^m\,e^{\pi i\,\s(\CM\, x\,b\,{}^tx)}f(x),\\
\left( {\widetilde R}_\CM\big(g_\CM(\alpha;t)\big)f\right)(x)&=& t^m\,\big( \det \alpha\big)^{\frac m2}\,f(x\,{}^t\alpha),
\end{eqnarray}
\begin{equation}
\left( {\widetilde R}_\CM\big(\s_{n,\CM;t}\big)f\right)(x)=t^m\,\left(
{\frac 1i}\right)^{\frac{mn}2} \big( \det \CM\big)^{\frac
n2}\,\int_{\rmn}f(y)\,e^{-2\,\pi i\,\s(\CM\,y\,{}^tx)}\,dy.
\end{equation}

We denote by $L^2_+\big(\rmn\big)$\,$\big(
\textrm{resp.}\,\,L^2_-\big(\rmn\big)\big)$ the subspace of
$L^2\big(\rmn\big)$ consisting of even (resp.\,odd) functions in
$L^2\big(\rmn\big)$. According to Formulas (3.10)--(3.12),
$R_{2,\CM}$ is decomposed into representations of $R_{2,\CM}^\pm$

\begin{equation*}
R_{2,\CM}=R_{2,\CM}^+\oplus R_{2,\CM}^-,
\end{equation*}
where $R_{2,\CM}^+$ and $R_{2,\CM}^-$ are the even Weil
representation and the odd Weil representation of $G$ that are
realized on $L^2_+\big(\rmn\big)$ and $L^2_-\big(\rmn\big)$
respectively. Obviously the center ${\mathscr Z}^J_{2,\CM}$ of
$G_{2,\CM}^J$ is given by
\begin{equation*}
{\mathscr Z}_{2,\CM}^J=\big\{ \big( (I_{2n},1),(0,0;\k)\big)\in
G_{2,\CM}^J\,\big\} \cong S(m,\BR).
\end{equation*}
We note that the restriction of $\omega_\CM$ to $G_{2,\CM}$
coincides with $R_{2,\CM}$ and $\omega_\CM(h)={\mathscr W}_\CM(h)$
for all $h\in \hrnm.$

\vskip 0.2cm\noindent
\begin{remark}
In the case $n=m=1,\
\omega_\CM$ is dealt in \cite{BS} and \cite{Ma}. We refer to
\cite{G} and \cite{KV} for more details about the Weil
representation $R_{2,\CM}$.
\end{remark}

\begin{remark}
The Schr{\"o}dinger-Weil representation is applied usefully to the theory of
Jacobi's sum \cite{Ma} and the theory of Maass-Jacobi forms \cite{P}.
\end{remark}

\end{section}

\vskip 1cm

\begin{section}
{{\bf Jacobi Forms of Half Integral Weight}}
\setcounter{equation}{0}

\newcommand\wg{\widetilde g}

\vskip 0.21cm

Let $\rho$ be a rational representation of
$GL(n,\mathbb{C})$ on a finite dimensional complex vector space
$V_{\rho}.$ Let ${\mathcal M}\in \mathbb R^{(m,m)}$ be a symmetric
half-integral semi-positive definite matrix of degree $m$. Let
$C^{\infty}(\BH_{n,m},V_{\rho})$ be the algebra of all
$C^{\infty}$ functions on $\BH_{n,m}$ with values in $V_{\rho}.$
For $f\in C^{\infty}(\BH_{n,m}, V_{\rho}),$ we define
\begin{align}
  & (f|_{\rho,{\mathcal M}}[(g,(\lambda,\mu;\kappa))])(\Om,Z) \notag \\
:= \,& e^{-2\,\pi\, i\,\sigma\left( {\mathcal M}(Z+\lambda\,
\Om+\mu)(C\Om+D)^{-1}C\,{}^t(Z+\lambda\,\Om\,+\,\mu)\right) }
\times e^{2\,\pi\, i\,\sigma\left( {\mathcal M}(\lambda\,
\Om\,{}^t\!\lambda\,+\,2\,\lambda\,{}^t\!Z+\,\kappa+
\mu\,{}^t\!\lambda) \right)} \\
&\times\rho(C\,\Om+D)^{-1}f(g\!\cdot\!\Om,(Z+\lambda\,
\Om+\mu)(C\,\Om+D)^{-1}),\notag
\end{align}
where $g=\left(\begin{matrix} A&B\\ C&D\end{matrix}\right)\in
Sp(n,\mathbb R),\ (\lambda,\mu;\kappa)\in H_{\mathbb R}^{(n,m)}$
and $(\Om,Z)\in \BH_{n,m}.$
\vspace{0.1in}\\
\noindent
\begin{definition}
Let $\rho$ and $\mathcal M$
be as above. Let
$$H_{\mathbb Z}^{(n,m)}:= \{ (\lambda,\mu;\kappa)\in H_{\mathbb R}^{(n,m)}\, \vert
\, \lambda,\mu\in \mathbb Z^{(m,n)},\ \kappa\in \mathbb
Z^{(m,m)}\,\ \}.$$ A $\textsf{Jacobi\ form}$ of index $\mathcal M$ with
respect to $\rho$ on a subgroup $\Gamma$ of $\Gamma_n$ of finite
index is a holomorphic function $f\in
C^{\infty}(\BH_{n,m},V_{\rho})$ satisfying the following
conditions (A) and (B):

\smallskip

\noindent (A) \,\ $f|_{\rho,{\mathcal M}}[\tilde{\gamma}] = f$ for
all $\tilde{\gamma}\in {\widetilde\Gamma}:= \Gamma \ltimes
H_{\mathbb Z}^{(n,m)}$.

\smallskip

\noindent (B) \,\ For each $M\in\Gamma_n$, $f|_{\rho,\CM}$ has a
Fourier expansion of the following form :
$$f(\Om,Z) = \sum\limits_{T=\,{}^tT\ge0\atop \text {half-integral}}
\sum\limits_{R\in \mathbb Z^{(n,m)}} c(T,R)\cdot e^{{ {2\pi
i}\over {\lambda_\G}}\,\sigma(T\Om)}\cdot e^{2\pi i\,\sigma(RZ)}$$

with a suitable $\lambda_\G\in\BZ^\times$ and
$c(T,R)\ne 0$ only if $\left(\begin{matrix} { 1\over {\lambda_\G}}T & \frac 12R\\
\frac 12\,^t\!R&{\mathcal M}\end{matrix}\right) \geqq 0$.
\end{definition}

\medskip

\indent If $n\geq 2,$ the condition (B) is superfluous by K{\"
o}cher principle\,(\,cf.\,\cite{Zi} Lemma 1.6). We denote by
$J_{\rho,\mathcal M}(\Gamma)$ the vector space of all Jacobi forms
of index $\mathcal{M}$ with respect to $\rho$ on $\Gamma$.
Ziegler\,(\,cf.\,\cite{Zi} Theorem 1.8 or \cite{EZ} Theorem 1.1\,)
proves that the vector space $J_{\rho,\mathcal {M}}(\Gamma)$ is
finite dimensional. In the special case $\rho(A)=(\det(A))^k$ with
$A\in GL(n,\BC)$ and a fixed $k\in\BZ$, we write $J_{k,\CM}(\G)$
instead of $J_{\rho,\CM}(\G)$ and call $k$ the {\it weight} of the
corresponding Jacobi forms. For more results on Jacobi forms with
$n>1$ and $m>1$, we refer to \cite{YJ6}-\cite{YJ9} and \cite{Zi}. Jacobi forms play an
important role in lifting elliptic cusp forms to Siegel cusp forms of degree $2n$ (cf.\,\cite{Ik}).

\vskip 0.2cm
\begin{definition}
A Jacobi form $f\in J_{\rho,\mathcal {M}}(\Gamma)$ is said to be
a $\textsf{cusp}$\,(\,or $\textsf{cuspidal\,) form}$ if $\begin{pmatrix} {
1\over {\lambda_\G}}T & {\frac 12}R\\ {\frac 12}\,^t\!R & \mathcal
{M}\end{pmatrix} > 0$ for any $T,\,R$ with $c(T,R)\ne 0.$ A Jacobi
form $f\in J_{\rho,\mathcal{M}}(\Gamma)$ is said to be
$\textsf{singular}$ if it admits a Fourier expansion such that a Fourier
coefficient $c(T,R)$ vanishes unless $\text{det}\begin{pmatrix} {
1\over {\lambda_\G}}T &{\frac 12}R\\ {\frac 12}\,^t\!R & \mathcal
{M}\end{pmatrix}=0.$
\end{definition}

\vskip 0.2cm
\begin{remark}
Singular Jacobi forms were characterized by a certain differential operator and the weight by the
author \cite{YJ8}.
\end{remark}

%%%%%%%%%%%%%%%%%%%%%%%%%%%%%%%%%%%%%%%%%%%%%%%%%%%%%%%%%%%%%%%%%%%%%%%%%%%%%%%%%%%%%%%%%%%%%%%%%%%%%%%%%%%%%%%%%%%%%%%%%%%%%%%%%%%%%%%%%%%%%%%%%%%%%%

\smallskip
Without loss of generality we may assume that $\rho$ is
irreducible. Then we choose a hermitian inner product $\langle\ ,\
\rangle$ on $V_{\rho}$ that is preserved under the unitary group
$U(n)\subset GL(n,\BC).$ For two Jacobi forms $f_1$ and $f_2$ in
$J_{\rho,\CM}(\Gamma)$, we define the Petersson inner product
formally by
\begin{equation}
(f_1,f_2):=\int_{\Gamma_{n,m}\backslash \BH_{n,m}}\langle\,
\rho(Y^{\frac 12})f_1(\Omega,Z),\rho(Y^{\frac
12})f_2(\Omega,Z)\rangle\,\k_\CM(\Omega,Z)\,dv.
\end{equation}
Here
\begin{equation}
dv=(\det Y)^{-(n+m+1)}[dX]\wedge [dY]\wedge [dU]\wedge [dV]
\end{equation}
is a $G^J$-invariant volume element on $\BH_{n,m}$
and
\begin{equation}
\kappa_\CM (\Om,Z):=e^{-4\pi i\,\sigma (\, {}^t({\rm Im}\,Z)\,\CM\, {\rm Im}\,Z\,({\rm Im}\,\Om)^{-1})}=
e^{-4\pi i\,\sigma (\, {}^tV\CM VY^{-1})},
\end{equation}
where $\Om=X+i\,Y\in \BH_n,\ Z=U+i\,V\in \BC^{(m,n)},\ X=(x_{ij}),\,Y=(y_{ij}),\,U=(u_{kl}),\,V=(v_{kl})$ real and
$$[dX]=\bigwedge_{i\leq j}dx_{ij},\quad  [dY]=\bigwedge_{i\leq j}dy_{ij},\quad
[dU]=\bigwedge_{k\leq l}du_{kl}\quad {\rm and}\quad[dV]=\bigwedge_{k\leq l}dv_{kl}.$$
A Jacobi form
$f$ in $J_{\rho,\CM}(\Gamma)$ is said to be {\it square\ integrable} if
$(f,f) < \infty.$ We note that cusp Jacobi forms are
square integrable and that $(f_1,f_2)$ is finite if
one of $f_1$ and $f_2$ is a cusp Jacobi form (cf.\,\cite{Zi},\
p.\,203).

\smallskip
For $g=\begin{pmatrix} A & B \\ C & D \end{pmatrix}\in G$, we set
\begin{equation}
J(g,\Om)=C\Om+D,\quad \Om\in\BH_n.
\end{equation}
\noindent
We define the map $J_\CM:G^J\times \BH_{n,m}\lrt \BC^\times$ by
\begin{equation}
J_\CM\big({\widetilde g},(\Om,Z)\big):=e^{2\pi i\,\sigma ( \CM [Z+\lambda \Om+\mu](C\Om+D)^{-1}C )}\cdot
e^{-2\pi i\,\sigma ( \CM(\lambda\,\Om\,{}^t\lambda\,+\,2\lambda\,{}^tZ+\kappa+\mu\,{}^t\lambda) )},
\end{equation}
where ${\widetilde g}=(g,(\lambda,\mu;\kappa))\in G^J$ with
$g=\begin{pmatrix} A & B \\ C & D \end{pmatrix}\in G$ and $(\lambda,\mu;\kappa)\in \hrnm.$ Here we use the
Siegel's notation $S[X]:=\,{}^tXSX$ for two matrices $S$ and $X$.

\vskip 0.21cm
We define the map $J_{\rho,\CM}:G^J\times\BH_{n,m}\lrt
GL(V_\rho)$ by
\begin{equation}
J_{\rho,\CM}(\wg,(\Om,Z))=J_\CM(\wg,(\Om,Z))\,\rho(J(g,\Omega)),
\end{equation}
where $\wg=(g,h)\in G^J$ with $g\in G$ and $h\in\hrnm.$ For a
function $f$ on $\BH_n$ with values in $V_\rho$, we can lift $f$
to a function $\Phi_f$ on $G^J$\,:
\begin{eqnarray*}
\Phi_f(\s):&=& (f|_{\rho,\CM}[\s])(iI_n,0)\\
&=& J_{\rho,\CM}(\s,(iI_n,0))^{-1}f(\s\!\cdot\! (iI_n,0)),\quad
\s\in G^J.
\end{eqnarray*}
A characterization of $\Phi_f$ for a cusp Jacobi form $f$ in
$J_{\rho,\CM}(\Gamma)$ was given by Takase
\cite[pp.\,162--164]{Ta1}.

%%%%%%%%%%%%%%%%%%%%%%%%%%%%%%%%%%%%%%%%%%%%%%%%%%%%%%%%%%%%%%%%%%%%%%%%%%%%%%%%%%%%%%%%%%%%%%%%%%%%%%%%%%%%%%%%%%%%%%%%%%%%%%%%%%%%%%%%%%%%%%%%%%%%%%%%%

\medskip We allow a weight $k$ to be half-integral.
Let
$$\mathfrak{S}=\left\{\, S\in\BC^{(n,n)}\,|\ S=\,{}^tS,\ {\rm Re}\,(S)>0\,\right\}$$
be a connected simply connected complex manifold. Then there is a uniquely determined holomorphic
function ${\rm det}^{1/2}$ on $\mathfrak{S}$ such that
\begin{eqnarray}
\left( {\rm det}^{1/2} S\right)^2 &=& \det S \qquad \textrm{for\ all}\ S\in \mathfrak{S},\\
{\rm det}^{1/2} S &=& (\det S)^{1/2} \qquad \textrm{for\ all}\ S\in \mathfrak{S}\cap \BR^{(n,n)}.
\end{eqnarray}
For each integer $k\in\BZ$ and $S\in \mathfrak{S}$, we put
$$ {\rm det}^{k/2} S=\left( {\rm det}^{1/2} S\right)^k.$$

\vskip 0.21cm
For brevity, we set $G=Sp(n,\BR).$
For any $g\in G$ and $\Om,\Om'\in\BH_n$, we put
\begin{eqnarray}
\varepsilon (g;\Om',\Om)&=& {\rm det}^{-{1/2}}\left( { {g\!\cdot\!\Om'-\overline{g\!\cdot\!\Om}}\over {2\,i}} \right)
{\rm det}^{1/2}\left( { {\Om'-\overline{\Om}}\over {2\,i}} \right) \\
&&\ \times |\det J(g,\Om')|^{-1/2}\,|\det J(g,\Om)|^{-1/2}.\nonumber
\end{eqnarray}

\vskip 0.21cm
For each $\Om\in\BH_n$, we define the function $\beta_\Om :G\times G\lrt T$ by
\begin{equation}
\beta_\Om (g_1,g_2)=\epsilon (g_1;\Om,g_2(\Om)),\quad g_1,g_2\in G.
\end{equation}
Then $\beta_\Om$ satisfies the cocycle condition and the cohomology class of $\beta_\Om$ of order two\,:
\begin{equation}
\beta_\Om (g_1,g_2)^2=\alpha_\Om (g_2)\, \alpha_\Om (g_1g_2)^{-1}\, \alpha_\Om (g_1),
\end{equation}
where
\begin{equation}
\alpha_\Om (g)={ {\det J(g,\Om)}\over {|\det J(g,\Om)|} }\,,\quad g\in G,\ \Om\in \BH_n.
\end{equation}
For any $\Om\in\BH_n,$ we let
\begin{equation*}
G_\Om=\left\{\,(g,\epsilon)\in G\times T\,|\ \epsilon^2=\alpha_\Om (g)^{-1}\,\right\}
\end{equation*}
be the two-fold covering group with the multiplication law
\begin{equation*}
(g_1,\epsilon_1)(g_2,\epsilon_2)=\big( g_1g_2,\,\epsilon_1 \epsilon_2 \beta_\Om (g_1,g_2) \big).
\end{equation*}
The covering group $G_\Om$ depends on the choice of $\Om\in\BH_n,$ i.e., the choice of a maximal compact
subgroup of $G$. However for any two element $\Om_1,\Om_2\in \BH_n,\ G_{\Om_1}$ is isomorphic to
$G_{\Om_2}$ (cf.\,\cite{Ta3}). We put
\begin{equation}
G_*:=G_{iI_n}.
\end{equation}
Takase \cite[p.\,131]{Ta3} defined the automorphic factor $J_{1/2}:G_*\times \BH_n\lrt \BC^\times$ by
\begin{equation}
J_{1/2}(g_{\epsilon},\Om):=\epsilon^{-1} \varepsilon (g;\Om, iI_n) |\det J(g,\Om)|^{1/2},
\end{equation}
where $g_{\epsilon}=(g,\epsilon)\in G_*$ with $g\in G$ and $\Om\in\BH_n.$ It is easily checked that
\begin{equation}
J_{1/2}(g_*h_*,\Om)=J_{1/2}(g_*,h\!\cdot\!\Om) J_{1/2}(h_*,\Om)
\end{equation}
for all $g_*=(g,\epsilon),\ h_*=(h,\eta)\in G_*$ and $\Om\in\BH_n$. Moreover
\begin{equation}
J_{1/2}(g_*,\Om)^2= \det (C\Om+D)
\end{equation}
for all $g_*=(g,\epsilon)\in G_*$ with $g=\begin{pmatrix} A & B \\ C & D \end{pmatrix}\in G.$

\vskip 0.21cm
Let $\pi_*:G_*\lrt G$ be the projection defined by $\pi_* (g,\epsilon)=g$. Let $\G$ be a subgroup of
the Siegel modular group $\G_n$ of finite index. Let $\G_*=\pi_*^{-1}(\G)\subset G_*.$ Let $\chi$ be a
finite order unitary character of $\G_*$. Let $k\in\BZ^+$ be a positive integer. We say that a holomorphic
function $\phi:\BH_n\lrt \BC$ is a $\textsf{Siegel\ modular\ form}$ of a $\textsf{half-integral\ weight}$ $k/2$ with
level $\G$ if it satisfies the condition
\begin{equation}
\phi (\g_*\!\cdot\Om)=\chi(\g_*) J_{1/2}(\g_*,\Om)^k\phi(\Om)
\end{equation}
for all $\g_*\in \G_*$ and $\Om\in\BH_n$. We denote by $M_{k/2}(\G,\chi)$ be the vector space of all
Siegel modular forms of weight $k/2$ with level $\G$. Let $S_{k/2}(\G,\chi)$ be the subspace of $M_{k/2}(\G,\chi)$
consisting of $\phi\in M_{k/2}(\G,\chi)$ such that
$$ |\phi(\Om)|\det ({\rm Im}\,\Om)^{k/4}\ {\rm is\ bounded\ on}\ \BH_n.$$
An element of $S_{k/2}(\G,\chi)$ is called a $\textsf{Siegel cusp form}$ of weight $k/2.$
The Petersson norm on $S_{k/2}(\G,\chi)$ is defined by
\begin{equation*}
||\phi||^2=\int_{\G\ba \BH_n}|\phi(\Om)|^2\,\det({\rm Im}\,\Om)^{k/2}\,dv_\Om,
\end{equation*}
where
$$dv_\Om =(\det Y)^{-(n+1)} [dX][dY]$$
is a $G$-invariant volume element on $\BH_n.$

\vskip 0.21cm
\begin{remark}
Using the Schr{\"o}dinger-Weil representation, Takase \cite{Ta4} established a bijective correspondence between the space of cuspidal
Jacobi forms and the space of Siegel cusp forms of half integral weight which is compatible with the action of Hecke operators.
For example, if $m$ is a positive integer, the classical result (cf.\,\cite{EZ} and \cite{Ib})
\begin{equation*}
J_{m,1}^{cusp}(\G_n)\cong S_{m-1/2}(\G_0(4))
\end{equation*}
can be obtained by the method of the representation theory. Here $\G_0(4)$ is the Hecke subgroup of the Siegel modular group
$\G_n$ and $J_{m,1}^{cusp}(\G_n)$ denotes the vector space of cuspidal Jacobi forms of weight $m$ and index $1$.
\end{remark}

\vskip 0.321cm
Now we are in a position to define the notion of Jacobi forms of half integral weight as follows.
\begin{definition}
Let $\G\subset \G_n$ be a subgroup of finite index. We put $\G_*=\pi_*^{-1}(\G)$ and
\begin{equation*}
\G_*^J=\G_*\ltimes H_\BZ^{(n,m)}.
\end{equation*}
\noindent
A holomorphic function $f:\BH_{n,m}\lrt
\BC$ is said to be a $\textsf{Jacobi form}$ of a $\textsf{weight}$ $k/2\in {\frac 12}\BZ$ (k\,:\,odd)
with $\textsf{level}$ $\G$ and $\textsf{index}$ $\CM$ for the character $\chi$ of $\Gamma_*^J$ of if it satisfies the following
transformation formula
\begin{equation}
f({\widetilde\g}_*\cdot
(\Om,Z))=\,\chi(\gamma_*)\,J_{k,\CM}({\widetilde\g}_*,(\Om,Z))f(\Om,Z)\quad
\textrm{for\ all}\ {\widetilde\g}_*\in \G^J_*.
\end{equation}

\noindent Here $J_{k,\CM}:\G^J_*\times\BH_{n,m}\lrt\BC$
is an automorphic factor defined by

\begin{eqnarray}
J_{k,\CM}\big( {\widetilde \g}_*,(\Om,Z)\big):&=& e^{2\,  \pi\,
i\,\sigma\big({\mathcal{M}}(Z+\lambda\,
\Om+\mu)(C\Om+D)^{-1}C\,{}^t(Z+\lambda \,\Om+\mu) \big)} \\
& &\  \times\,\, e^{-2\pi i\sigma\left( {\mathcal{M}}(\lambda\,
\Om\,{}^t\!\lambda\,+\,2\lambda\,{}^t\!Z\,+\,\kappa\,+\,
\mu\,{}^t\!\lambda)\right) }
J_{1/2}(\gamma_*,\Omega)^k,\nonumber
%\hskip 1cm
\end{eqnarray}
where $\widetilde{\gamma}_*=(\gamma_*,(\lambda,\mu;\kappa))\in \G^J_*$ with $\gamma_*=(\gamma,\epsilon)\in \G_*,$
$\gamma=\begin{pmatrix} A&B\\
C&D\end{pmatrix}\in \G,\  (\lambda,\mu,\kappa)\in H_{\BZ}^{(n,m)}$
and $(\Om,Z)\in \BH_{n,m}.$

\end{definition}

%\begin{equation}
%A
%\end{equation}

\end{section}

\vskip 1cm

\newcommand\mfm{{\mathscr F}^{(\CM)} }
\newcommand\mfoz{{\mathscr F}^{(\CM)}_{\Om,Z} }
\newcommand\wg{{\widetilde g} }
\newcommand\wgm{{\widetilde \gamma} }
\newcommand\Tm{\Theta^{(\CM)} }

\begin{section}
{{\bf Covariant Maps for the Schr{\"o}dinger-Weil representation }}
\setcounter{equation}{0}

\vskip 0.2cm As before we let $\CM$ be a symmetric positive
definite $m\times m$ real matrix. We keep the notations in the previous sections.

\vskip 0.21cm
We define the mapping ${\mathscr
F}^{(\CM)}:\BH_{n,m}\lrt L^2\big(\rmn\big)$ by

\begin{equation}
\mfm (\Om,Z)(x)=\,e^{\pi i\,\s\{ \CM
(x\,\Om\,{}^tx+\,2\,x\,{}^tZ)\} },\quad (\Om,Z)\in\BH_{n,m},\
x\in\rmn.
\end{equation}

\noindent For brevity we put $\mfoz:=\mfm (\Om,Z)$ for $(\Om,Z)\in
\BH_{n,m}.$ Takase \cite{Ta3} proved that $G_{2,\CM}$ is isomorphic to $G_{iI_n}$ (cf.\,(3.6) and (4.14)). Therefore
we will use $G_*:=G_{iI_n}$ instead of $G_{2,\CM}.$

\vskip 0.21cm
We set
\begin{equation}
G_*^J:= G_*\ltimes H_\BR^{(n,m)}.
\end{equation}
We note that $G_*^J$ acts on $\BH_{n,m}$ via the canonical projection of $G_*^J$ onto $G^J.$

\vskip 0.21cm
Now we assume that $m$ is odd.
We define the automorphic factor $J_\CM^*:G^J_*\times \BH_{n,m}\lrt
\BC^{\times}$ for $G^J_*$ on $\BH_{n,m}$ by

\begin{eqnarray}
   J_\CM^*({\widetilde g}_*,(\Om,Z))&= e^{\pi
   i\,\sigma\left({\mathcal{M}}(Z+\lambda\,
\Om+\mu)(C\Om+D)^{-1}C\,{}^t(Z+\l\,\Om+\mu)\right)} \hskip 3cm \\
& \times e^{-\pi i\,\sigma \left( \mathcal{M}(\lambda\,
\Om\,^t\!\lambda\,+\,2\,\lambda\,{}^t\!Z\,+\,\kappa\,+\,
\mu\,{}^t\!\lambda ) \right) } J_{1/2}((g,\epsilon),\Om)^m,\nonumber
\end{eqnarray}
where ${\widetilde g}_*=((g,\epsilon),(\lambda,\mu;\kappa))\in G^J_*$ with $g=\begin{pmatrix} A&B\\
C&D\end{pmatrix}\in G,\ (\lambda,\mu;\kappa)\in
H_{\BR}^{(n,m)}$ and $(\Om,Z)\in \BH_{n,m}.$

\begin{theorem} Let $m$ be an odd positive integer. Then the map ${\mathscr F}^{(\CM)}:\BH_{n,m}\lrt
L^2\big(\rmn\big)$ defined by (5.1) is a $\textsf{covariant map}$ for the
Schr{\"o}dinger-Weil representation $\omega_\CM$ of $G^J_*$ and the
automorphic factor $J_\CM^*$ for $G^J_*$ on $\BH_{n,m}$ defined by
Formula (5.3). In other words, $\mfm$ satisfies the following
covariance relation

\begin{equation}
\om_\CM ({\widetilde g}_*)\mfoz=J_\CM^* \big( {\widetilde g}_*,(\Om,Z)\big)^{-1}
\mfm_{{\widetilde g_*}\cdot (\Om,Z)}
\end{equation}

\noindent for all ${\widetilde g}_*\in G^J_*$ and $(\Om,Z)\in
\BH_{n,m}.$
\end{theorem}

\noindent {\it Proof.} For an element
$\wg_*=((g,\epsilon),(\lambda,\mu;\kappa))\in G^J_*$ with $g=\begin{pmatrix} A &
B
\\ C & D \end{pmatrix}\in Sp(n,\BR),$ we put
$(\Om_*,Z_*)=\wg_*\cdot (\Om,Z)$ for $(\Om,Z)\in\BH_{n,m}.$ Then we
have
\begin{eqnarray*}
&\Om_*=g\cdot \Om=(A\Om+B)(C\Om+D)^{-1},\\
&Z_*=(Z+\l\, \Om+\mu)(C\Om+D)^{-1}.
\end{eqnarray*}

\newcommand\ep{\epsilon}

\noindent In this section we use the notations $t(b),\
g(\alpha)$ and $\sigma_n$ in Section 3. Since the following
elements $h(\l,\mu;\kappa)_\ep,\ t(b;\ep),\ g(\alpha,\ep)$ and $\s_{n,\ep}$ of $G^J_*$
defined by

\begin{eqnarray*}
 h(\l,\mu;\kappa)_\ep&=&\big( (I_{2n},\ep),(\l,\mu;\kappa))\quad \textrm{with}\
\ep=\pm \,1,\ \l,\mu\in\rmn,\ \kappa\in\BR^{(m,m)},\\
t(b;\ep)&=&\big( (t(b),\ep),(0,0;0)\big)\quad \textrm{with}\ \ep=\pm \,1,\ b\in S(n,\BR),\\
 g(\alpha;\ep)&=&\big( (g(\alpha),\ep),(0,0;0)\big)\quad \textrm{with}\
%\alpha\in GL(n,\BR),\
\ep=\pm\, 1\ {\rm if}\ \det \alpha >0,\ \ep=\pm \,i \ {\rm if}\ \det \alpha <0, \\
\s_{n,\ep}&=&\big( (\s_n,\ep),(0,0;0)\big)\ {\rm with}\ \ep^2=(-i)^n.
\end{eqnarray*}

\noindent generate the group $G^J_*$, it suffices to prove the
covariance relation (5.4) for the above generators.

\vskip 0.5cm\noindent {\bf Case I.} $\wg_*=h(\l,\mu;\kappa)_\ep$ with
$\ep=\pm\,1$ and $\l,\mu\in\rmn,\ \kappa\in\BR^{(m,m)}.$

\vskip 0.1cm In this case, we have
$$\Om_*=\Om,\quad Z_*=Z+\la\,\Om+\mu.$$
We put
\begin{equation*}
h(\l,\mu;\kappa)_+:=\big( (I_{2n},1),(\l,\mu;\kappa)\big)
\end{equation*}
and
\begin{equation*}
h(\l,\mu;\kappa)_-:=\big( (I_{2n},-1),(\l,\mu;\kappa)\big).
\end{equation*}
It is easily seen that according to Formula (4.15), we obtain
\begin{equation*}
J_{1/2}((I_{2n},1),\Om)=1 \quad {\rm and}\quad J_{1/2}((I_{2n},-1),\Om)=-1.
\end{equation*}
Therefore we get
\begin{equation*}
J_\CM^*\big( h(\l,\mu;\kappa)_+,(\Om,Z)\big)=\,e^{-\pi i\,\s\{ \CM(\l\,\Om\,{}^t\l
+2\,\l\,{}^tZ+\kappa+\mu\,{}^t\!\l)\} }
\end{equation*}
and
\begin{equation*}
J_\CM^*\big( h(\l,\mu;\kappa)_-,(\Om,Z)\big)=\,-e^{-\pi i\,\s\{ \CM(\l\,\Om\,{}^t\l
+2\,\l\,{}^tZ+\kappa+\mu\,{}^t\!\l)\} }.
\end{equation*}

According to Formulas (2.4), (3.8) and (3.9), for $x\in \rmn,$

\begin{eqnarray*}
& & \left( \om_\CM \big( h(\l,\mu;\kappa)_+\big) \mfoz\right) (x)\\
&=& e^{\pi i\,\sigma\{ \CM(\kappa+\mu\,{}^t\l+2\,x\,{}^t\mu)\} }
\mfoz
(x+\l)\\
&=& e^{\pi i\,\sigma\{ \CM(\kappa+\mu\,{}^t\l+2\,x\,{}^t\mu)\}
}\,e^{\pi i \,\s\{ \CM(
(x+\l)\Om\,{}^t(x+\l)+\,2\,(x+\l)\,{}^tZ)\} }.
\end{eqnarray*}

\noindent On the other hand, according to Formula (5.3), for $x\in
\rmn,$

\begin{eqnarray*}
& & J_\CM^* \big( h(\l,\mu;\kappa)_+,(\Om,Z)\big)^{-1}\mfm_{\wg_*\cdot
(\Om,Z)}(x)\\
&=& J_\CM^* \big(
h(\l,\mu;\kappa)_+,(\Om,Z)\big)^{-1}\mfm_{\Om,Z+\l\,\Om+\mu}(x)\\
&=& e^{\pi i\,\sigma\{ \CM(\l
\Om\,{}^t\!\l\,+\,2\,\l\,{}^tZ\,+\,\kappa\,+\,\mu\,{}^t\!\l)\}
}\cdot e^{\pi
i\,\sigma\{ \CM (x\, \Om\,{}^tx\,+\,2\,x\,{}^t(Z+\l\Om+\mu))\} }\\
&=& e^{\pi i\,\sigma\{ \CM(\kappa+\mu\,{}^t\l+2\,x\,{}^t\mu)\}
}\,e^{\pi i \,\s\{ \CM(
(x+\l)\Om\,{}^t(x+\l)+\,2\,(x+\l)\,{}^tZ)\} }.
\end{eqnarray*}

\noindent Therefore we prove the covariance relation (5.4) in the
case $\wg_*=h(\l,\mu;\kappa)_+$.

\vskip 0.21cm
Similarly we can prove the covariance relation (5.4) in the $\wg_*=h(\l,\mu;\kappa)_-.$ In fact,
\begin{eqnarray*}
& & \left( \om_\CM \big( h(\l,\mu;\kappa)_-\big) \mfoz\right) (x)\\
&=&\,- e^{\pi i\,\sigma\{ \CM(\kappa+\mu\,{}^t\l+2\,x\,{}^t\mu)\}
}\,e^{\pi i \,\s\{ \CM(
(x+\l)\Om\,{}^t(x+\l)+\,2\,(x+\l)\,{}^tZ)\} }\\
&=& J_\CM^* \big( h(\l,\mu;\kappa)_-,(\Om,Z)\big)^{-1}\mfm_{\wg_*\cdot
(\Om,Z)}(x).
\end{eqnarray*}

\vskip 0.5cm\noindent {\bf Case II.} $\wg_*=t(b;\ep)$ with $\ep=\pm\,1$ and
$b\in S(n,\BR).$

\vskip 0.1cm In this case, we have
$$\Om_*=\Om+b \quad {\rm and}\quad Z_*=Z.$$
We put
$$ t(b)_+=\big( (t(b),1),(0,0;0)\big)$$
and
$$ t(b)_-=\big( (t(b),-1),(0,0;0)\big).$$
It is easily seen that
\begin{equation*}
J_{1/2}((t(b),1),\Om)=1\quad {\rm and}\quad J_{1/2}((t(b),-1),\Om)=-1.
\end{equation*}
Therefore we get
\begin{equation*}
J_\CM^*\big( t(b)_+,(\Om,Z) \big)=1 \quad {\rm and}\quad J_\CM^*\big( t(b)_-,(\Om,Z) \big)=-1.
\end{equation*}

\noindent According to Formulas (3.8) and (3.10), we obtain
\begin{equation*}
\left(\om_\CM \big( t(b)_+\big) \mfoz \right) (x)=\,e^{\pi\,i\,\sigma (\CM\,xb\,{}^tx)} \mfoz (x),\quad
x\in\rmn.
\end{equation*}

\noindent On the other hand, according to Formula (5.3), for $x\in
\rmn,$ we obtain

\begin{eqnarray*}
& & J_\CM^* \big( t(b)_+,(\Om,Z)\big)^{-1}\mfm_{t(b)_+\cdot
(\Om,Z)}(x)\\
&=& \mfm_{\Om+b,Z}(x)\\
&=& e^{\pi i\,\sigma\left( \CM \left(
x(\Om+b)\,{}^tx+2\,x\,{}^tZ\right)\right)}\\
&=& \,e^{\pi\,i\,\sigma (\CM\,xb\,{}^tx)} \mfoz (x).
\end{eqnarray*}

\noindent Therefore we prove the covariance relation (5.4) in the
case $\wg_*=t(b)_+$ with $b\in S(n,\BR).$

\vskip 0.21cm
Similarly we can prove the covariance relation (5.4) in the $\wg_*=t(b)_-.$ In fact,

\begin{eqnarray*}
& & \left(\om_\CM \big( t(b)_-\big) \mfoz \right) (x)\\
&=&\,-e^{\pi\,i\,\sigma (\CM\,xb\,{}^tx)} \mfoz (x)\\
&=& J_\CM^* \big( t(b)_-,(\Om,Z)\big)^{-1}\mfm_{t(b)_-\cdot
(\Om,Z)}(x).
\end{eqnarray*}

\vskip 0.52cm\noindent {\bf Case III.} \ \ $\wg_*=\big( (g(\alpha),\ep),(0,0;0)\big)$ with
$\ep=\pm\,1\ (\det\, \alpha>0)$ and $\ep=\pm\, i\ (\det\,\alpha < 0).$ \\
\indent \ \ \ \ \ \  \indent\indent Here $\alpha\in GL(n,\BR).$

\vskip 0.2311cm In this case, we have
$$\Om_*=\,{}^t\alpha\,\Om\,\alpha\quad {\rm and} \quad Z_*=Z\alpha.$$
We put
\begin{eqnarray*}
 g(\alpha)_+:&=& \big( (g(\alpha),1),(0,0;0)\big),\\
 g(\alpha)_-:&=& \big( (g(\alpha),-1),(0,0;0)\big),\\
 g(\alpha)^+:&=& \big( (g(\alpha),i),(0,0;0)\big),\\
 g(\alpha)^-:&=& \big( (g(\alpha),-i),(0,0;0)\big).
\end{eqnarray*}
And we can show easily that
\begin{eqnarray*}
J_{1/2}\big( (g(\alpha),1),\Om)&=&(\det\,\alpha)^{-1/2},\\
J_{1/2}\big( (g(\alpha),-1),\Om)&=&-(\det\,\alpha)^{-1/2},\\
J_{1/2}\big( (g(\alpha),i),\Om)&=& i\,(\det\,\alpha)^{-1/2},\\
J_{1/2}\big( (g(\alpha),-i),\Om)&=& -i\,(\det\,\alpha)^{-1/2}.
\end{eqnarray*}

Using Formulas (3.5), (3.11) and (5.3), we can show the covariance relation (5.4) in the
case $\wg_*=\big( (g(\alpha),\ep),(0,0;0)\big)$ with $\ep=\pm \,1,\ \pm\,i$ and $\alpha\in GL(n,\alpha).$

%\vskip 3cm
%\noindent According to Formula (4.7), for $x\in \rmn,$
%\begin{eqnarray*}
%& & \left( \om_\CM \big(\wg \big) \mfoz\right) (x)\\
%&=& (\det\alpha)^{{\frac m2}} \mfoz
%(x\,{}^t\alpha)\\
%&=&\,(\det\alpha)^{{\frac m2}}\cdot e^{\pi i\,\sigma\{
%\CM(x\,{}^t\alpha\,\Om\,{}^t(x\,{}^t\alpha)+2\,x\,{}^t\alpha\,{}^tZ)\}}.
%\end{eqnarray*}

%\noindent On the other hand, according to Formula (5.2), for $x\in \rmn,$
%\begin{eqnarray*}
%& & J_\CM \big( \wg,(\Om,Z)\big)^{-1}\mfm_{\wg\cdot (\Om,Z)}(x)\\
%&=&(\det\alpha)^{{\frac m2}}\mfm_{{}^t\alpha\,\Om\,\alpha,Z\alpha}(x)\\
%&=&\,(\det\alpha)^{{\frac m2}}\cdot e^{\pi i\,\sigma\{\CM(x\,{}^t\alpha\,\Om\,{}^t(x\,{}^t\alpha)+2\,x\,{}^t\alpha\,{}^tZ)\}}.
%\end{eqnarray*}

\vskip 0.52cm\noindent
{\bf Case IV.} $\wg_*=\big( (\sigma_n,\ep),(0,0;0)\big)$ with
$\sigma_n=\begin{pmatrix} 0 & -I_n \\ I_n & \ 0 \end{pmatrix}$ and $\ep^2=(-i)^n.$

\vskip 0.21cm In this case, we have
$$\Om_*=-\Om^{-1}\quad {\rm and} \quad Z_*=Z\,\Om^{-1}.$$

In order to prove the covariance relation (5.3), we need
the following useful lemma.

\begin{lemma} For a fixed element $\Om\in \BH_n$ and a fixed
element $Z\in\BC^{(m,n)},$ we obtain the following property
\begin{equation}
\int_{\rmn} e^{\pi\,i\,\sigma
(x\,\Om\,{}^tx+2\,x\,{}^tZ)}dx_{11}\cdots dx_{mn} = \left( \det
{\Omega\over i}\right)^{-{\frac m2}}\,
e^{-\pi\,i\,\sigma(Z\,\Om^{-1}\,{}^tZ)},
\end{equation}
where $x=(x_{ij})\in \BR^{(m,n)}.$
\end{lemma}

\noindent {\it Proof of Lemma 5.1.} By a simple computation, we
see that
$$e^{\pi i\, \sigma ( x\Om\, {}^tx +
2x\,{}^tZ )}= e^{-\pi i\,\sigma (Z\Om^{-1}\,{}^tZ )}\cdot e^{\pi
i\,\sigma \{(x+Z\Om^{-1})\Om\,{}^t(x+Z\Om^{-1})\} }.$$ Since the
real Jacobi group $Sp(n,\BR)\ltimes H_\BR^{(m,n)}$ acts on
${\mathbb H}_{n,m}$ holomorphically, we may put
$$\Om=\,i\,A\,{}^t\!A,\quad Z=iV,\quad\  A\in\BR^{(n,n)},\quad
V=(v_{ij})\in\BR^{(m,n)}.$$
Then we obtain
\begin{eqnarray*}
& & \int_{\BR^{(m,n)}}  e^{\pi i\, \sigma ( x\Om\, {}^tx +
2x\,{}^tZ )} dx_{11}\cdots dx_{mn} \\
&=& e^{-\pi i\,\sigma (Z\Omega^{-1}\,{}^tZ)} \int_{\BR^{(m,n)}}
e^{\pi i\,\sigma [\{
x+iV(iA\,{}^t\!A)^{-1}\}(iA\,{}^t\!A)\,{}^t\!\{
x+iV(iA\,{}^t\!A)^{-1}\} ]}\,dx_{11}\cdots dx_{mn}\\
&=&e^{-\pi i\,\sigma (Z\Omega^{-1}\,{}^tZ)} \int_{\BR^{(m,n)}}
e^{\pi i\,\sigma [\{ x+V(A\,{}^t\!A)^{-1}\}A\,{}^t\!A\,{}^t\!\{
x+V(A\,{}^t\!A)^{-1}\} ]}\,dx_{11}\cdots dx_{mn}\\
&=& e^{-\pi i\,\sigma (Z\Omega^{-1}\,{}^tZ)} \int_{\BR^{(m,n)}}
e^{-\pi \,\sigma\{ (uA)\,{}^t\!(uA)\} }\,du_{11}\cdots
du_{mn}\quad \big(\,{\rm Put}\ u=
x+V(A\,{}^t\!A)^{-1}=(u_{ij}) \,\big)\\
&=& e^{-\pi i\,\sigma (Z\Omega^{-1}\,{}^tZ)} \int_{\BR^{(m,n)}}
e^{-\pi \,\sigma (w\,{}^t\!w)} (\det A)^{-m}\,dw_{11}\cdots
dw_{mn}\quad \big(\,{\rm Put}\ w=uA=(w_{ij})\,\big)\\
&=& e^{-\pi i\,\sigma (Z\Omega^{-1}\,{}^tZ)} \, (\det A)^{-m}\cdot
\left( \prod_{i=1}^m \prod_{j=1}^g \int_\BR e^{-\pi\,
w_{ij}^2}\,dw_{ij}\right)\\
&=& e^{-\pi i\,\sigma (Z\Omega^{-1}\,{}^tZ)} \, (\det A)^{-m}\quad
\big(\,{\rm because}\ \int_\BR e^{-\pi\,
w_{ij}^2}\,dw_{ij}=1\quad {\rm for\ all}\ i,j\,\big)\\
&=& e^{-\pi i\,\sigma (Z\Omega^{-1}\,{}^tZ)} \, \left( \det \big(
A\, {}^t\!A \big)\right)^{-{\frac m2}}\\
&=& e^{-\pi i\,\sigma (Z\Omega^{-1}\,{}^tZ)} \, \left( \det \left(
{ {\Omega}\over i } \right)\right)^{-{\frac m2}}.
\end{eqnarray*}

\noindent This completes the proof of Lemma 5.1. \hfill $\square$

\vskip 0.21cm According to Formulas (3.8) and (3.12), for $x\in\rmn,$ we obtain

\begin{eqnarray*}
& & \left( \om_\CM \big(\wg_* \big) \mfoz\right) (x)\\
&=& \ep^m\,\left( {\frac 1i}\right)^{{mn}\over 2}\big( \det
\CM\big)^{\frac n2}\,\int_{\rmn} \mfoz (y)\,e^{-2\pi\,i\,\s\,(\CM
\,y\,{}^tx)}dy\\
&=& \ep^m\,\left( {\frac 1i}\right)^{{mn}\over 2}\big( \det
\CM\big)^{\frac n2}\,\int_{\rmn} e^{\pi\,i\,\sigma \{
\CM(y\,\Om\,{}^ty+2\, y\,{}^tZ)\}
}\,e^{-2\pi\,i\,\sigma(\CM \,y\,{}^tx)}dy\\
&=&\ep^m\,\left( {\frac 1i}\right)^{{mn}\over 2}\big( \det
\CM\big)^{\frac n2}\,\int_{\rmn} e^{ \pi \,i\,\sigma\left\{
\CM\left( y\,\Om\,{}^ty\,+\,2\,y\,{}^t(Z-x)\right)\right\} }dy.
\end{eqnarray*}

\noindent If we substitute $u=\CM^{\frac 12}\,y,$ then $du=\left(
\det \CM\right)^{\frac n2}\,dy.$ Therefore according to Lemma 5.1,
we obtain

\begin{eqnarray*}
& & \left( \om_\CM \big(\wg_* \big) \mfoz\right) (x)\\
&=& \ep^m\,\left( {\frac 1i}\right)^{{mn}\over 2}\big( \det
\CM\big)^{\frac n2}\,\int_{\rmn} e^{\pi\,i\,\sigma \left(
u\,\Om\,{}^tu\,+\,2\,\CM^{1/2}\,u\,{}^t(Z-x)\right)}\,\left(\det \CM\right)^{-{\frac n2}}du\\
&=&\ep^m\,\left( {\frac 1i}\right)^{{mn}\over 2}\,\int_{\rmn}
e^{\pi\,i\,\sigma \left(
u\,\Om\,{}^tu\,+\,2\,u\,{}^t(\CM^{1/2}\,(Z-x))\right)}\,du\\
&=&\ep^m\,\left( {\frac 1i}\right)^{{mn}\over 2}\, \left( \det
{\Omega\over i}\right)^{-{\frac m2}}\, e^{-\pi\,i\,\sigma\left\{
\CM^{1/2} (Z-x)\,\Om^{-1}\,{}^t(Z-x)\,\CM^{1/2}\right\}}\quad (\textrm{by\ Lemma\ 5.1})\\
&=&\ep^m\,\,\left(\det \Om\right)^{-{\frac m2}}\,
e^{-\pi\,i\,\sigma\left( \CM\,(Z-x)\,\Om^{-1}\,{}^t(Z-x)\right)
}\\
&=&\ep^m\,\left(\det \Om\right)^{-{\frac m2}}\,
e^{-\pi\,i\,\sigma\left(
\CM(Z\,\Om^{-1}\,{}^tZ\,+\,x\,\Om^{-1}\,{}^tx\,-\,2\,Z\,\Om^{-1}\,{}^tx)\right) }.
\end{eqnarray*}

\noindent On the other hand,
\begin{eqnarray*}
& & J_{1/2}\big( (\sigma_n,\ep),\Om\big)\\
&=& \ep^{-1}\,{\det}^{-1/2}\left( { {\s_n\!\cdot\!\Om- \overline{\s_n\!\cdot\!(iI_n)}}\over {2\,i} }\right)
\,{\det}^{-1/2}\left( { {\Om- \overline{(iI_n)}}\over {2\,i} }\right)\,
|J(\s_n,iI_n)|^{-1/2}\\
&=& \ep^{-1}\,{\det}^{-1/2}\left( { {-i\,\Om^{-1}\big(\Om-\overline{i\,I_n}\big)}\over {2\,i} }\right)
\,{\det}^{-1/2}\left( { {\Om- \overline{(iI_n)}}\over {2\,i} }\right)\,(i^n)^{-1/2}\\
&=& \ep^{-1}\,{\det}^{-1/2}\left( -i\,\Om^{-1}\right)\,i^{-n/2} \\
&=& \ep^{-1}\,\left( {\frac 1i}\right)^{-n/2}\,(\det\,\Om)^{1/2}\,\,i^{-n/2}\\
&=& \ep^{-1}\,(\det \,\Om)^{1/2}.
\end{eqnarray*}

\noindent
Therefore, according to Formula (5.3), for $x\in \rmn,$ we obtain

\begin{eqnarray*}
& & J_\CM^* \big( \wg_*,(\Om,Z)\big)^{-1}\mfm_{\wg_*\cdot
(\Om,Z)}(x)\\
&=&\,e^{-\pi\,i\,\sigma(\CM\,Z\,\Om^{-1}\,{}^tZ)}\,J_{1/2}\big( (\s_n,\ep),\Om\big)^{-m}\, \mfm_{-\Om^{-1},\,Z\,\Om^{-1}}(x)\\
&=&\,\ep^m\,\left(\det \Om\right)^{-{\frac m2}}\,e^{-\pi\,i\,\sigma(\CM\,Z\,\Om^{-1}\,{}^tZ)}\,
e^{\pi\,i\,\sigma\left\{ \CM \left(
x\,(-\Om^{-1})\,{}^tx\,+\,2\,x\,{}^t(Z\,\Om^{-1}) \right)\right\}} \\
&=&\,\ep^m\,\left(\det \Om\right)^{-{\frac m2}}\,
e^{-\pi\,i\,\sigma\left(
\CM(Z\,\Om^{-1}\,{}^tZ\,+\,x\,\Om^{-1}\,{}^tx\,-\,2\,Z\,\Om^{-1}\,{}^tx)\right)
}.
\end{eqnarray*}

\noindent Hence we prove the covariance relation (5.4) in the
case $\wg_*=\big( (\sigma_n,\ep),(0,0;0)\big)$ with $\ep^2=(-i)^n.$  Since $J_\CM^*$ is an automorphic factor for
$G^J_*$ on $\BH_{n,m}$, we see that if the covariance relation (5.4)
holds for two elements $\wg_*,\,{\widetilde h}_*$ in $G^J_*$, then it holds
for $\wg_*{\widetilde h}_*.$ Finally we complete the proof. \hfill $\square$

\vskip 0.37cm
It is natural to raise the following question\,:
\vskip 0.23cm\noindent
{\bf Problem\,:} Find all the covariant maps for the Sch{\"o}dinger-Weil representation $\omega_\CM$ on
$\BH_{n,m}$ and the automorphic factor $J_\CM^*.$

\end{section}

\vskip 1cm

\begin{section}{{\bf Construction of Jacobi Forms of Half Integral Weight }}
\setcounter{equation}{0}

\vskip 0.3cm Let $(\pi,V_\pi)$ be a unitary representation of
$G^J_*$ on the representation space $V_\pi$. Let $\G$ be an arithmetic subgroup of the Siegel modular group $\G_n.$
We set $\G_*=\pi_*^{-1}(\G)$ and
$$\G_*^J=\G_* \ltimes H_\BZ^{(n,m)}.$$
We assume that
$(\pi,V_\pi)$ satisfies the following conditions (A) and (B):
\vskip 0.21cm \noindent {\bf (A)} There exists a vector valued map
\begin{equation*}
{\mathscr F}:\BH_{n,m}\lrt V_\pi,\quad\ (\Om,Z)\mapsto {\mathscr
F}_{\Om,Z}:={\mathscr F}(\Om,Z)
\end{equation*}

\noindent satisfying the following covariance relation
\begin{equation}
\pi\big( {\widetilde g}_*\big) {\mathscr F}_{\Om,Z}=\psi\big(
{\widetilde g}_*\big)\,J_*\big( {\widetilde g}_*,(\Om,Z)\big)^{-1}\,{\mathscr F}_{ {\widetilde g}_*\cdot (\Om,Z)
}\quad \textrm{for all}\ {\widetilde g}_*\in G^J_*,\
(\Om,Z)\in\BH_{n,m},
\end{equation}

\noindent where $\psi$ is a character of $G^J_*$ and $J_*:G^J_*\times
\BH_{n,m}\lrt \BC^\times$ is a certain automorphic factor for $G^J_*$
on $\BH_{n,m}.$

\vskip 0.1cm \noindent {\bf (B)} There exists a linear functional
$\theta:V_\pi\lrt \BC$ which is $\textsf{ semi-invariant}$ under the action of
$\G^J_*$, in other words, for all ${\widetilde\g}_*\in
\G^J_*$ and $(\Om,Z)\in\BH_{n,m},$

\begin{equation}
\langle\, \pi^*\big( {\widetilde\g}_*\big)\theta,\,{\mathscr F}_{\Om,Z}\,\rangle=\langle\,
\theta,\pi\big( {\widetilde\g}_*\big)^{-1}{\mathscr F}_{\Om,Z}\,\rangle=\chi\big({\widetilde\g}_*\big)\,\langle\,
\theta,\,{\mathscr F}_{\Om,Z}\,\rangle ,
\end{equation}

\noindent where $\pi^*$ is the contragredient of $\pi$ and
$\chi:\G^J_*\lrt T$ is a unitary character of
$\G^J_*$.

\vskip 0.2cm Under the assumptions (A) and (B) on a unitary
representation $(\pi,V_\pi)$, we define the function $\Theta$ on
$\BH_{n,m}$ by

\begin{equation}
\Theta(\Om,Z):=\,\langle\,\theta,{\mathscr
F}_{\Om,Z}\,\rangle=\theta\big({\mathscr F}_{\Om,Z}\big),\quad\
(\Om,Z)\in\BH_{n,m}.
\end{equation}

We now shall see that $\Theta$ is an automorphic form on
$\BH_{n,m}$ with respect to $\G^J_*$ for the automorphic
factor $J_*$.

\begin{lemma} Let $(\pi,V_\pi)$ be a unitary representation of
$G^J_*$ satisfying the above assumptions (A) and (B). Then the
function $\Theta$ on $\BH_{n,m}$ defined by (6.3) satisfies the
following modular transformation behavior
\begin{equation}
\Theta\big({\widetilde \g}_*\cdot(\Om,Z)\big)=\,\psi\big({\widetilde\g}_*\big)^{-1}
\,\chi\big({\widetilde\g}_*\big)^{-1}\,J_*\big( {\widetilde \g}_*,(\Om,Z)\big)\,\Theta(\Om,Z)
\end{equation}

\noindent for all ${\widetilde\g}_*\in \G_*^J$ and
$(\Om,Z)\in\BH_{n,m}.$
\end{lemma}

\noindent {\it Proof.} For any ${\widetilde\g}_*\in \G_*^J$
and $(\Om,Z)\in\BH_{n,m},$ according to the assumptions (6.1) and
(6.2), we obtain
\begin{eqnarray*}
& & \Theta\big({\widetilde\g}_*\cdot(\Om,Z)\big)=\big\langle\,\theta,
{\mathscr F}_{ {\widetilde \g}_*\cdot (\Om,Z)} \big\rangle \hskip 5cm\\
&=&\big\langle\,\theta,\psi\big( {\widetilde \g}_*\big)^{-1}J_*\big(
{\widetilde \g}_*,(\Om,Z)\big)\,\pi\big({\widetilde\g}_*\big){\mathscr F}_{ \Om,Z}\,\big\rangle\\
&=&\psi\big( {\widetilde \g}_*\big)^{-1} J_*\big( {\widetilde\g}_*,(\Om,Z)\big)\,\big\langle
\,\theta,\pi\big({\widetilde\g}_*\big){\mathscr F}_{ \Om,Z}\,\big\rangle\\
&=&\,\psi\big( {\widetilde\g}_*\big)^{-1}\,\chi\big({\widetilde\g}_*\big)^{-1}
\,J_*\big( {\widetilde \g}_*,(\Om,Z)\big)\,\big\langle\,\theta,{\mathscr F}_{ \Om,Z}\,\big\rangle\\
&=&\,\psi\big( {\widetilde\g}_*\big)^{-1}\,\chi\big({\widetilde\g}_*\big)^{-1}
\,J_*\big( {\widetilde \g}_*,(\Om,Z)\big)\,\Theta(\Om,Z).
\end{eqnarray*}
\hfill$\square$

\newcommand\zmn{\BZ^{(m,n)} }
\newcommand\wgam{\widetilde\gamma}

\vskip 0.2cm Now for a positive definite integral symmetric matrix
$\CM$ of degree $m$, we define the holomorphic function
$\Theta_\CM:\BH_{n,m}\lrt\BC$ by

\begin{equation}
\Theta_\CM (\Om,Z):=\sum_{\xi\in \BZ^{(m,n)}}
e^{\pi\,i\,\sigma\left( \CM(
\xi\,\Om\,{}^t\xi\,+\,2\,\xi\,{}^tZ)\right) },\quad (\Om,Z)\in
\BH_{n,m}.
\end{equation}

\begin{theorem}
Let $m$ be an odd positive integer.
Let $\CM$ be a symmetric positive definite integral matrix of degree $m$ such that $\det\,\CM=1$.
Let $\G$ be an arithmetic subgroup of $\G_n$ generated by all the following elements
\begin{equation*}
t(b)=\begin{pmatrix} I_n & b \\ 0 & I_n \end{pmatrix},\quad
g(\alpha)=\begin{pmatrix} {}^t\alpha & 0 \\ 0 & \alpha^{-1} \end{pmatrix},\quad
\s_n=\begin{pmatrix} 0 & -I_n \\ I_n & \ 0\end{pmatrix},
\end{equation*}
where $b=\,{}^tb\in S(n,\BZ)$ with even diagonal and $\alpha\in GL(n,\BZ).$
Then for any $\wgam_*\in\Gamma^J_*$, the function
$\Theta_\CM$ satisfies the functional equation
\begin{equation}
\Theta_\CM\big( \wgam_*\!\cdot\! (\Om,Z)\big)=\rho_\CM (\wgam_*)\,
J_\CM^*\big( \wgam_*,(\Om,Z)\big) \Theta_\CM(\Om,Z),\quad (\Om,Z)\in
\BH_{n,m},
\end{equation}
\end{theorem}
\noindent where $\rho_\CM$ is a character of $\G^J_*$
and $J_\CM^*:G^J_*\times \BH_{n,m}\lrt \BC^{\times}$ is the automorphic factor for $G^J_*$ on $\BH_{n,m}$
defined by the formula (5.3).

\newcommand{\ep}{\epsilon}

\noindent {\it Proof.}  For an element
$\wgam_*=\big( (\g,\ep),(\lambda,\mu;\kappa)\big)\in \G^J_*$ with $\g=\begin{pmatrix}
A & B \\ C & D \end{pmatrix}\in \G$ and $(\l,\mu;\kappa)\in H_\BZ^{(n,m)},$ we put
$(\Om_*,Z_*)=\wgam_*\!\cdot\! (\Om,Z)$ for $(\Om,Z)\in\BH_{n,m}.$ Then
we have
\begin{eqnarray*}
&\Om_*=\g\cdot \Om=(A\Om+B)(C\Om+D)^{-1},\\
&Z_*=(Z+\l\, \Om+\mu)(C\Om+D)^{-1}.
\end{eqnarray*}

\noindent We define the linear functional $\vartheta$ on $L^2\big(
\rmn\big)$ by

\begin{equation*}
\vartheta (f)=\langle \vartheta,f \rangle:=\sum_{\xi\in
\zmn}f(\xi),\quad\ f\in L^2\big( \rmn\big).
\end{equation*}

\noindent We note that $\Theta_\CM(\Om,Z)=\vartheta\big(
\mfoz\big).$ Since $\mfm$ is a covariant map for the
Schr{\"o}dinger-Weil representation $\omega_\CM$ by Theorem 5.1,
according to Lemma 6.1, it suffices to prove that $\vartheta$ is
semi-invariant for $\om_\CM$ under the action of $\G^J_*$, in other
words, $\vartheta$ satisfies the following semi-invariance
relation
\begin{equation}
\Big\langle\, \vartheta,\om_\CM\big( \wgam_*\big)\mfoz \,\Big\rangle
=\,\rho_\CM \big( \wgam_*\big)^{-1} \,\Big\langle\, \vartheta,\mfoz
\,\Big\rangle
\end{equation}

\noindent for all $\wgam_*\in \G^J_*$ and $(\Om,Z)\in \BH_{n,m}.$

\vskip 0.212cm We note that the following elements
$h(\lambda,\mu;\kappa)_\ep,\ t(b;\ep),\,g(\alpha;\ep)$ and $\s_{n,\ep}$ of $\G^J_*$
defined by
\begin{eqnarray*}
&& h(\la,\mu;\kappa)_\ep=\big( (I_{2n},\ep),(\la,\mu;\kappa)\big)\
\textrm{with}\ \ep=\pm\,1,\ \la,\mu\in \BZ^{(m,n)}\ \textrm{and}\ \kappa\in\BZ^{(m,m)} ,\\
&&t(b;\ep)=\big( (t(b),\ep),(0,0;0)\big)\ \textrm{with}\ \ep=\pm\,1\ \textrm{and}
                   \ b=\,{}^tb\in S(n,\BZ)\ {\rm even\ diagonal},\\
&& g(\alpha;\ep)=\big( (g(\alpha),\ep),(0,0;0)\big)\
\textrm{with}\ \ep=\pm\,1,\ \pm\,i\ \textrm{and}\ \alpha\in GL(n,\BZ),\\
 &&\s_{n,\ep}=\big( (\s_n,\ep),(0,0;0)\big)\ {\rm with}\ \ep^2=(-i)^n.
\end{eqnarray*}
generate the Jacobi modular group $\G^J_*.$ Therefore it suffices to
prove the semi-invariance relation (6.7) for the above generators
of $\G^J_*.$

\vskip 0.5cm\noindent {\bf Case I.} $\wgam_*=h(\l,\mu;\kappa)_\ep$ with
$\ep=\pm\,1,\ \l,\mu\in\zmn,\ \kappa\in\BZ^{(m,m)}.$

\vskip 0.1cm In this case, we have
$$\Om_*=\Om\quad {\rm and} \quad Z_*=Z+\la\,\Om+\mu.$$
We put
\begin{equation*}
h(\l,\mu;\kappa)_+:=\big( (I_{2n},1),(\l,\mu;\kappa)\big)
\end{equation*}
and
\begin{equation*}
h(\l,\mu;\kappa)_-:=\big( (I_{2n},-1),(\l,\mu;\kappa)\big).
\end{equation*}
It is easily seen that according to Formula (4.15), we obtain
\begin{equation*}
J_{1/2}((I_{2n},1),\Om)=1 \quad {\rm and}\quad J_{1/2}((I_{2n},-1),\Om)=-1.
\end{equation*}
Therefore we get
\begin{equation*}
J_\CM^*\big( h(\l,\mu;\kappa)_+,(\Om,Z)\big)=\,e^{-\pi i\,\s\{ \CM(\l\,\Om\,{}^t\l
+2\,\l\,{}^tZ+\kappa+\mu\,{}^t\!\l)\} }
\end{equation*}
and
\begin{equation*}
J_\CM^*\big( h(\l,\mu;\kappa)_-,(\Om,Z)\big)=\,-e^{-\pi i\,\s\{ \CM(\l\,\Om\,{}^t\l
+2\,\l\,{}^tZ+\kappa+\mu\,{}^t\!\l)\} }.
\end{equation*}
According to the covariance relation (5.4),
\begin{eqnarray*}
& &\big\langle \,\vartheta, \om_\CM\big( h(\l,\mu;\kappa)_+\big)
\mfoz\,\big\rangle\\
&=&\,\big\langle \,\vartheta,
J_\CM^*\big( h(\l,\mu;\kappa)_+,(\Om,Z)\big)^{-1}\mfm_{h(\l,\mu;\kappa)_+\cdot(\Om,Z)}\,\big\rangle\\
&=&\,J_\CM^* \big( h(\l,\mu;\kappa)_+,(\Om,Z)\big)^{-1}\,\big\langle \,\vartheta,
\mfm_{\Om,Z+\la\,\Om+\mu}\,\big\rangle\\
&=&\,J_\CM^* \big( h(\l,\mu;\kappa)_+,(\Om,Z)\big)^{-1}\sum_{A\in\zmn}e^{\pi\,i\,\s\left\{ \CM
\left( A\Om\,{}^t\!A\,+\,2\,A\,{}^t(Z+\l\,\Om+\mu)\right) \right\} }\\
&=&\,J_\CM^* \big( h(\l,\mu;\kappa)_+,(\Om,Z)\big)^{-1}\cdot
e^{-\pi\,i\,\sigma\left( \CM
(\l\,\Om\,{}^t\!\l\,+\,2\,\l\,{}^tZ)\right)}\\
& &\times \sum_{A\in\zmn} e^{2\,\pi\,i\,\sigma(\CM A\,{}^t\mu)}
e^{\pi\,i\,\s\left\{ \CM \left(
(A+\l)\,\Om\,{}^t\!(A+\l)\,+\,2\,(A+\l)\,{}^tZ\right) \right\}
}\\
&=& e^{\pi\,i\,\sigma\left( \CM(\kappa\,+\,\mu\,{}^t\!\l)\right)}
\,\big\langle \,\vartheta, \mfm_{\Om,Z}\,\big\rangle.
\end{eqnarray*}

\noindent Here we used the fact that $\s(\CM A\,{}^t\mu)$ is an
integer. In a similar way, we get
\begin{equation*}
\big\langle \,\vartheta, \om_\CM\big( h(\l,\mu;\kappa)_-\big)
\mfoz\,\big\rangle=-e^{\pi\,i\,\sigma\left( \CM(\kappa\,+\,\mu\,{}^t\!\l)\right)}
\,\big\langle \,\vartheta, \mfm_{\Om,Z}\,\big\rangle.
\end{equation*}
We put
\begin{equation*}
\rho_\CM\big( h(\l,\mu;\kappa)_+\big)= e^{-\pi\,i\,\sigma\left( \CM(\kappa\,+\,\mu\,{}^t\!\l)\right)}
\end{equation*}
and
\begin{equation*}
\rho_\CM\big( h(\l,\mu;\kappa)_-\big)= -e^{-\pi\,i\,\sigma\left( \CM(\kappa\,+\,\mu\,{}^t\!\l)\right)}.
\end{equation*}
Therefore $\vartheta$
satisfies the semi-invariance relation (6.7) in the case $\wgam_*=h(\l,\mu;\kappa)_\ep$ with
$\ep=\pm\,1,\ \l,\mu\in\zmn$ and $\kappa\in\BZ^{(m,m)}.$

\vskip 0.5cm\noindent {\bf Case II.} $\wgam_*=t(b;\ep)$ with
$\ep=\pm\,1,\ b=\,{}^tb\in S(n,\BZ)$ even diagonal.

\vskip 0.1cm In this case, we have
$$\Om_*=\Om+b\quad {\rm and}\quad Z_*=Z.$$
We put
$$ t(b)_+=\big( (t(b),1),(0,0;0)\big)$$
and
$$ t(b)_-=\big( (t(b),-1),(0,0;0)\big).$$
It is easily seen that
\begin{equation*}
J_{1/2}((t(b),1),\Om)=1\quad {\rm and}\quad J_{1/2}((t(b),-1),\Om)=-1.
\end{equation*}
Therefore we get
\begin{equation*}
J_\CM^*\big( t(b)_+,(\Om,Z) \big)=1 \quad {\rm and}\quad J_\CM^*\big( t(b)_-,(\Om,Z) \big)=-1.
\end{equation*}

\noindent According to the covariance relation (5.4), we obtain
\begin{eqnarray*}
& &\big\langle \,\vartheta, \om_\CM\big( t(b)_+\big)
\mfoz\,\big\rangle\\
&=&\,\big\langle \,\vartheta,
J_\CM^*\big( t(b)_+,(\Om,Z)\big)^{-1}\mfm_{t(b)_+\cdot(\Om,Z)}\,\big\rangle\\
&=&\,\big\langle \,\vartheta,
\mfm_{\Om+b,Z}\,\big\rangle\\
&=&\, \sum_{A\in\zmn} e^{\pi\,i\,\s\left\{ \CM \left(
A\,(\Om+b)\,{}^t\!A\,+\,2\,A\,{}^tZ\right) \right\}
}\\
&=&\, \sum_{A\in\zmn} e^{\pi\,i\,\s\left( \CM \left(
A\,\Om\,{}^t\!A\,+\,2\,A\,{}^tZ\right) \right) }\cdot
e^{\pi\,i\,\sigma(\CM A\,b\,{}^t\!A)}\\
&=&\, \sum_{A\in\zmn} e^{\pi\,i\,\sigma\left( \CM (
A\,\Om\,{}^t\!A\,+\,2\,A\,{}^tZ ) \right) }\\
&=&\,\big\langle \,\vartheta, \mfm_{\Om,Z}\,\big\rangle.
\end{eqnarray*}

\noindent Here we used the fact that $\s(\CM A\,b\,{}^t\!A)$ is an
even integer because $b$ is even integral. In a similar way, we obtain
\begin{equation*}
\big\langle \,\vartheta, \om_\CM\big( t(b)_-\big)
\mfoz\,\big\rangle=-\big\langle \,\vartheta, \mfm_{\Om,Z}\,\big\rangle.
\end{equation*}
We put
\begin{equation*}
\rho_\CM \big( t(b)_+ \big)=1 \quad {\rm and}\quad
\rho_\CM \big( t(b)_- \big)=-1.
\end{equation*}
Therefore $\vartheta$ satisfies the semi-invariance
relation (6.7) in the case $\wgam_*=t(b;\ep)$ with
$\ep=\pm\,1,\ b\in S(n,\BZ).$

\vskip 0.52cm\noindent {\bf Case III.} \ \ $\wg_*=\big( (g(\alpha),\ep),(0,0;0)\big)$ with
$\ep=\pm\,1\ (\det\, \alpha>0)$ and $\ep=\pm\, i\ (\det\,\alpha < 0).$ \\
\indent \ \ \ \ \ \  \indent\indent Here $\alpha\in GL(n,\BZ).$

\vskip 0.2311cm In this case, we have
$$\Om_*=\,{}^t\alpha\,\Om\,\alpha\quad {\rm and} \quad Z_*=Z\alpha.$$
We put
\begin{eqnarray*}
 g(\alpha)_+:&=& \big( (g(\alpha),1),(0,0;0)\big),\\
 g(\alpha)_-:&=& \big( (g(\alpha),-1),(0,0;0)\big),\\
 g(\alpha)^+:&=& \big( (g(\alpha),i),(0,0;0)\big),\\
 g(\alpha)^-:&=& \big( (g(\alpha),-i),(0,0;0)\big).
\end{eqnarray*}
And we can show easily that
\begin{eqnarray*}
J_{1/2}\big( (g(\alpha),1),\Om)&=&(\det\,\alpha)^{-1/2},\\
J_{1/2}\big( (g(\alpha),-1),\Om)&=&-(\det\,\alpha)^{-1/2},\\
J_{1/2}\big( (g(\alpha),i),\Om)&=& i\,(\det\,\alpha)^{-1/2},\\
J_{1/2}\big( (g(\alpha),-i),\Om)&=& -i\,(\det\,\alpha)^{-1/2}.
\end{eqnarray*}
\noindent According to the covariance relation (5.4), we obtain
\begin{eqnarray*}
& &\big\langle \,\vartheta, \om_\CM\big( g(\alpha)_+\big)
\mfoz\,\big\rangle\\
&=&\,\big\langle \,\vartheta,
J_\CM^*\big( g(\alpha)_+,(\Om,Z)\big)^{-1}\mfm_{g(\alpha)_+\cdot(\Om,Z)}\,\big\rangle\\
&=&\,\left( \det \alpha \right)^{\frac m2}\,\big\langle
\,\vartheta,
\mfm_{{}^t\alpha\,\Om\,\alpha,Z\,\alpha}\,\big\rangle\\
&=&\,\left( \det \alpha \right)^{\frac m2} \sum_{A\in\BZ^{(m,n)} }
\mfm_{{}^t\alpha\,\Om\,\alpha,Z\,\alpha}(A)\\
&=&\,\left( \det \alpha \right)^{\frac m2}\sum_{A\in\BZ^{(m,n)} }
e^{\pi i\,\sigma\{
\CM(A\,{}^t\alpha\,\Om\,{}^t(A\,{}^t\alpha)+2\,A\,{}^t\alpha\,{}^tZ)\}} \\
&=&\,\left( \det \alpha \right)^{\frac m2}\,\big\langle
\,\vartheta, \mfm_{\Om,Z}\,\big\rangle.
\end{eqnarray*}
In a similar way we get
\begin{eqnarray*}
\big\langle \,\vartheta, \om_\CM\big( g(\alpha)_-\big)
\mfoz\,\big\rangle &=& (-1)^m\,\left( \det \alpha \right)^{\frac m2}\,\big\langle
\,\vartheta, \mfm_{\Om,Z}\,\big\rangle,\\
\big\langle \,\vartheta, \om_\CM\big( g(\alpha)^+\big)
\mfoz\,\big\rangle &=& i^m\,\left( \det \alpha \right)^{\frac m2}\,\big\langle
\,\vartheta, \mfm_{\Om,Z}\,\big\rangle,\\
\big\langle \,\vartheta, \om_\CM\big( g(\alpha)^-\big)
\mfoz\,\big\rangle &=& (-i)^m\,\left( \det \alpha \right)^{\frac m2}\,\big\langle
\,\vartheta, \mfm_{\Om,Z}\,\big\rangle.\\
\end{eqnarray*}
Now we put
\begin{eqnarray*}
\rho_\CM \big( g(\alpha)_+\big)&=& (\det\,\alpha )^{\frac m2},\\
\rho_\CM \big( g(\alpha)_-\big)&=& (-1)^m(\det\,\alpha )^{\frac m2},\\
\rho_\CM \big( g(\alpha)^+\big)&=& i^m\,(\det\,\alpha )^{\frac m2},\\
\rho_\CM \big( g(\alpha)^-\big)&=& (-i)^m\,(\det\,\alpha )^{\frac m2}.
\end{eqnarray*}

\noindent Therefore $\vartheta$
satisfies the semi-invariance relation (6.7) in the case
$\wgam_*=g(\alpha;\ep)$ with $\ep=\pm\,1,\ \pm\,i,\ \alpha\in GL(n,\BZ).$

\vskip 0.52cm\noindent
{\bf Case IV.} $\wgam_*=\big( (\sigma_n,\ep),(0,0;0)\big)$ with
$\sigma_n=\begin{pmatrix} 0 & -I_n \\ I_n & \ 0 \end{pmatrix}$ and $\ep^2=(-i)^n.$

\vskip 0.21cm In this case, we have
$$\Om_*=-\Om^{-1}\quad {\rm and} \quad Z_*=Z\,\Om^{-1}.$$

\noindent In the process of the proof of Theorem 5.1, using Lemma
5.1, we already showed that

\begin{equation}
\int_{\rmn} e^{\pi\,i\,\sigma( \CM(y\,\Om\,{}^ty\,+\,2\,y\,{}^tZ)
) } dy =\,\big(\det \CM\big)^{-{\frac n2}}\left( \det {\Om \over
i}\right)^{-{\frac
m2}}\,e^{-\pi\,i\,\sigma(\CM\,Z\,\Om^{-1}\,{}^tZ)}.
\end{equation}

\noindent By (6.8), we see that the Fourier transform of $\mfoz$
is given by

\begin{equation}
\widehat {\mfoz} (x)=\,\big(\det \CM\big)^{-{\frac n2}}\left( \det
{\Om \over i}\right)^{-{\frac m2}}\,e^{-\pi\,i\,\sigma(\CM\,(Z-x)\,\Om^{-1}\,{}^t(Z-x))},
\end{equation}
where $\widehat{f}$ is the Fourier transform of $f$ defined by
\begin{equation*}
\widehat{f} (x)=\int_\rmn f(y)\,e^{-2\,\pi\,i\,\s (y\,{}^tx)}\,dy,\quad x\in\rmn.
\end{equation*}
On the other hand, in the process of the proof of Case IV in Theorem 5.1, we showed that
\begin{equation*}
J_\CM^* \big( \widetilde{g}_*,(\Om,Z)\big)=\ep^{-m}\,(\det\,\Om)^{\frac m2}\,e^{\pi\,i\,\s (\CM Z\,\Om^{-1}{}^tZ)}.
\end{equation*}

\vskip 0.21cm
\noindent According to the covariance relation (5.4), Formula
(6.9) and Poisson summation formula, we obtain

\begin{eqnarray*}
& &\big\langle \,\vartheta, \om_\CM\big( \wgam_*\big)
\mfoz\,\big\rangle\\
&=&\,\big\langle \,\vartheta,
J_\CM^* \big( \wgam_*,(\Om,Z)\big)^{-1}\mfm_{\wgam_*\cdot(\Om,Z)}\,\big\rangle\\
&=&\,J_\CM^*\big( \wgam_*,(\Om,Z)\big)^{-1} \big\langle \,\vartheta,
\mfm_{-\Om^{-1},Z\Om^{-1}}\,\big\rangle\\
&=&\,\ep^m\,(\det \Om)^{-{\frac
m2}}\,e^{-\pi\,i\,\sigma(\CM\,Z\,\Om^{-1}\,{}^tZ) }
\sum_{A\in\BZ^{(m,n)}} e^{-\pi\,i\,\sigma\left( \CM(
A\,\Om^{-1}\,{}^t\!A\,-\,2\,A\,\Om^{-1}\,{}^tZ) \right)}\\
&=&\,\ep^m\,(\det \Om)^{-{\frac m2}}\, \sum_{A\in\BZ^{(m,n)}}
e^{-\pi\,i\,\sigma\left( \CM(Z\,\Om^{-1}\,{}^tZ\,+\,
A\,\Om^{-1}\,{}^t\!A\,-\,2\,A\,\Om^{-1}\,{}^tZ) \right)}\\
&=&\,\ep^m\,(\det \Om)^{-{\frac m2}}\, \sum_{A\in\BZ^{(m,n)}}
e^{-\pi\,i\,\sigma\left( \CM(Z-A)\,\Om^{-1}\,{}^t(Z-A)\right)}\\
&=&\,\ep^m\,(\det \Om)^{-{\frac m2}}\big(\det \CM\big)^{{\frac n2}}\left(
\det {\Om \over i}\right)^{{\frac m2}}\,\sum_{A\in\BZ^{(m,n)}}
\widehat {\mfoz} (\CM\,A) \quad (\,\textrm{by\ Formula}\ (6.9))\\
&=&\,\ep^m\,\left( \det {{I_n} \over
i}\right)^{{\frac m2}}\,\sum_{A\in\BZ^{(m,n)}}
\widehat{\mfoz} (A) \quad (\,\textrm{because}\ \det\,\CM=1) \\
&=&\,\ep^m\,\left( \det {{I_n} \over
i}\right)^{{\frac m2}}\,\sum_{A\in\BZ^{(m,n)}}
{\mfoz} (A) \quad (\,\textrm{by\ Poisson summation formula}) \\
&=&\,\ep^m\,(-i)^{{mn}\over 2}\,\big\langle \,\vartheta,
\mfm_{\Om,Z}\,\big\rangle.
\end{eqnarray*}

\noindent
We put
\begin{equation*}
\rho_\CM (\wgam_*)=\,\ep^m\,(-i)^{{mn}\over 2}.
\end{equation*}

\noindent Therefore $\vartheta$ satisfies the
semi-invariance relation (6.7) in the case $\wgam_*=\big( (\sigma_n,\ep),(0,0;0)\big)$
with $\ep^2=(-i)^n.$ The proof of Case IV is completed.

\vskip 0.21cm
Since $J_\CM^*$ is an automorphic
factor for $G^J_*$ on $\BH_{n,m}$, we see that if the formula (6.6)
holds for two elements $\wgam_1,\wgam_2$ in $\G^J$, then it holds
for $\wgam_1 \wgam_2.$ Finally we complete the proof of Theorem
6.1. \hfill $\square$

\begin{corollary}
Let $\G_*^J$ and $\rho_\CM$ be as before in Theorem 6.1. If $m$ is odd, then
$\Theta_\CM (\Om,Z)$ is a Jacobi form of a half integral weight ${\frac m2}$ and index ${{\CM}\over 2}$
with respect to an arithmetic subgroup $\G_*^J$ for a character $\rho_\CM$ of $\G_*^J$.
\end{corollary}

\vskip 0.21cm\noindent
\begin{remark}
Let $a=(a_1,a_2)\in\BZ^n\times \BZ^n$ with $a_1,a_2\in\BZ^n$. Takase \cite{Ta3} considered the following theta series
defined by
\begin{equation*}
\vartheta_a (\Om,Z):=\sum_{\ell\in\BZ^n}\, e^{\pi\,i\,((\ell+a_1)\,\Om\,{}^t(\ell+a_2)+\,2\,(\ell+a_1)\,{}^t(Z+a_2))},
\end{equation*}
where $\Om\in \BH_n$ and $Z\in \BC^n.$ We put
\begin{equation*}
\vartheta_a^* (\Om,Z):=\,e^{-a_1\,{}^ta_2}\,\vartheta_a (\Om,Z).
\end{equation*}
We let $\G_0$ be an arithmetic subgroup of $\G_n$ consisting of
$\g=\begin{pmatrix} A & B \\ C & D \end{pmatrix}$ such that
\vskip 0.127cm
(1) $L\g =L,$ where $L=\BZ^n\times \BZ^n.$
\vskip 0.12cm
(2) $(xA+y\,C)\,{}^t(xB+yD)\equiv x\,{}^y\ ({\rm mod}\ 2\BZ)\ {\rm for\ all}\ (x,y)\in L.$

\vskip 0.21cm
\noindent
We put
\begin{equation*}
\G_{*,0}=\pi_*^{-1}(\G_0)\quad {\rm and}\quad \G_{*,0}^J=\G_{*,0}\ltimes H_\BZ^{(n,1)}.
\end{equation*}
He proved that for any $\wgam_*=\big( (\g,\ep),(\l,\mu;t))\in \G_{*,0}^J$ with $\g\in \G_0,$
the following transformation formula
\begin{equation*}
\vartheta_{a\g^{-1}}^* \big( \wgam_*\!\cdot\! (\Om,Z)\big)=\rho( (\g,\ep))\,\chi_a((\l,\mu;t))\,
J_{\bf 1}(\wgam_*,(\Om,Z))\,\vartheta_a^* (\Om,Z)
\end{equation*}
holds, where $\chi_a$ denotes the unitary character of $L\times \BR$ defined by
\begin{equation*}
\chi_a ((\l,\mu;t))=\,e^{2\,\pi\,i\,\left( t\,+{\frac 12}\,\l\,{}^t\mu\,-\,\l\,{}^ta_2\,+\,a_1\,{}^t\mu\right)},\quad
(\l,\mu;t)\in L\times\BR
\end{equation*}
and ${\bf 1}=(1)$ denotes the $1\times 1$ matrix. Here $\rho:\G_{*,0}\lrt T$ is a certain unitary character
that is given explicitly in \cite[Theorem 5.3,\, p.134]{Ta3}.
\end{remark}

\end{section}

\vskip 1cm
\bibliography{central}

\end{document}